\documentclass[leqno]{article}
\usepackage[top=1.5in,bottom=1in,left=1.1in,right=1.1in,a4paper]{geometry}
\usepackage{amsfonts,latexsym,amsmath,amssymb,amsthm}
\usepackage{float,xcolor,array,fancyhdr,subfigure}
\usepackage[]{hyperref} 
\usepackage{enumitem,mathtools}
\usepackage{accents}
\usepackage{tikz,pgfplots,siunitx}
\usepackage{multirow}
\usepackage{algorithmic,algorithm2e,setspace}
\usepgfplotslibrary{groupplots}
\usetikzlibrary{positioning,arrows,fit}
\pgfplotsset{compat=newest}
\usepackage[textsize=small,colorinlistoftodos]{todonotes}
\newtheorem*{remark*}    {Remark}
\newtheorem {remark}     {Remark}
\newtheorem {theorem}    {Theorem}
\newtheorem {definition} {Definition}

\numberwithin{equation}{section}

\newcommand{\dd}{\mathrm{d}}

\newcommand{\ii}        {\mathrm{i}}

\newcommand{\misfit}    {\mathcal{J}}
\newcommand{\misfitI}   {\mathcal{I}}

\newcommand{\lagrangian}{\mathcal{L}}

\newcommand{\bx}        {x}

\newcommand{\Ree}       {\operatorname{Re}}
\newcommand{\Gobs}      {G_{obs}}
\newcommand{\Gc}        {G_{c}}
\newcommand{\Rc}        {R_{c}}
\newcommand{\Su}        {S_{\mathcal{U}_0}}
\newcommand{\Subar}     {\overline{\Su}}
\newcommand{\dc}        {\delta c}
\newcommand{\dGc}       {D_c\Gc}
\newcommand{\adjoint}{\gamma}
\newcommand{\tildadj}{\tilde{\adjoint}}

\newcommand*\samethanks[1][\value{footnote}]{\footnotemark[#1]}
\newcommand{\myblue}  {blue!70!white}
\newcommand{\mygreen} {green!50!black}
\newcommand{\myred}   {red!60!white}

\title{
Inverse problem for the Helmholtz equation with Cauchy data: 
reconstruction with conditional well-posedness driven iterative 
regularization
}
\author{{Giovanni Alessandrini\thanks{Dipartimento di Matematica e
    Geoscienze, Universit\`{a} di Trieste,
    Italy.}}
\and {Maarten V. de Hoop\thanks{Department of Computational and
    Applied Mathematics and Department of Earth Science, Rice
    University, Houston TX 77005, USA.}}
\and{Florian Faucher\thanks{Project-Team Magique-3D, 
     Inria Bordeaux Sud-Ouest Research Center,
     Laboratoire de Math\'ematiques et de leurs Applications, 
     Universit\'e de Pau et des Pays de l'Adour, UMR CNRS 5142, France.}}
\and {Romina Gaburro\thanks{
      Department of Mathematics and Statistics, 
      Health Research Institute (HRI), University of Limerick, Limerick, Ireland.}}
\and {Eva Sincich\samethanks[1]
}
}

\newlength {\modelwidth}    
 
\newcommand{\modelfile}{}   
  
\begin{document}

\maketitle

\begin{abstract}
In this paper, we study the performance of Full Waveform Inversion
(FWI) from time-harmonic Cauchy data via conditional well-posedness
driven iterative regularization. The Cauchy data can be obtained with
dual sensors measuring the pressure and the normal velocity. 
We define a novel misfit functional which, adapted to the Cauchy data,
allows the independent location of experimental and computational sources.
The conditional well-posedness is obtained for a hierarchy of subspaces in
which the inverse problem with partial data is Lipschitz stable. Here,
these subspaces yield piecewise linear representations of the wave
speed on given domain partitions. Domain partitions can be adaptively
obtained through segmentation of the gradient. The domain partitions
can be taken as a coarsening of an unstructured tetrahedral mesh
associated with a finite element discretization of the Helmholtz
equation. We illustrate the effectiveness of the iterative
regularization through computational experiments with
data in dimension three. In comparison with earlier work, the Cauchy
data do not suffer from eigenfrequencies in the configurations.
\end{abstract}

\section{Introduction}

Iterative methods for the recovery of subsurface parameters have been
collectively referred to, in reflection seismology, as Full Waveform
Inversion (FWI). FWI consists in the minimization of the residuals,
defined as some difference between the observed and modelled data. FWI
was originally formulated with time-domain data using an energy norm
in the misfit functional by Lailly \cite{Lailly1983} and
Tarantola \cite{Tarantola1984,Tarantola1987a}. 
The time-harmonic formulation of the seismic inverse problem was later
considered by Pratt \textit{et al.} \cite{Pratt1990,Pratt1996}. Standardly, 
one applies the adjoint state method for the implementation of the 
minimization procedure.

To mitigate the nonlinearity and ill-posedness of the inverse
problem, hierarchical multiscale strategies have been developed.  In
the time-domain, Bunks \textit{et al.} \cite{Bunks1995} proposed successive inversion of
data subsets of increasing frequency contents. This multiscale
approach can be related to the subspace search method introduced in
\cite{Kennett1988}. Multiscale Gauss-Newton-Krylov methods were
developed by Akcelik \textit{et al.} \cite{Akcelik2002} and
many further developments have taken place since then. 
The application of wavelet bases enabling successive
levels of model compression were considered by Loris \textit{et al.}
\cite{Loris2007,Loris2010} in wave-equation tomography and in FWI by
Lin \textit{et al.} in \cite{Lin2012}, Yuan and Simons in \cite{Yuan2014,Yuan2015}. 
These ideas were natural but lacked foundation and understanding of
convergence. In \cite{deHoop2012}, the authors developed an iterative 
regularization approach with a multilevel strategy derived from 
conditional Lipschitz stability estimates with a convergence analysis. For the
regularization they introduced a multiscale hierarchy of subspaces for
which the Lipschitz stability estimates hold, with an associated model
compression rate through projections. This compression rate appears in
the convergence analysis and mitigates the growth of stability
constants with scale refinement, that is, increasing dimension of the
subspaces. In the time-harmonic case, such quantitative stability
estimates have been proved for piecewise constant \cite{Beretta2012_AcousticStability}
and piecewise linear \cite{AdHGS} representations of the
wave speed with a given domain partition; 
we refer to \cite{AdHGS} for more extended bibliography.
Note that, in view of recent work by C{\^a}rstea \textit{et al.} 
\cite{Carstea2016}, 
on a germane problem in elasticity, it seems quite possible to 
recover the partition as well, in the case where the domain is 
composed of sub-analytic sets.
In the piecewise constant
case, this includes subspaces defined by Haar wavelets, where the analysis in \cite{AdHGS}
is adapted to general domain partitions such as unstructured
tetrahedral meshes, that can be associated with a segmentation. 
Here, we present a computational framework of our reconstruction algorithm
via iterative regularization using piecewise linear representations 
as stable subspaces, as well as experiments.

We use Cauchy data assimilated from dual sensor acquisition. 
Cauchy data do not suffer from eigenfrequencies unlike the
Dirichlet-to-Neumann map, which, in fact, cannot be observed directly
in seismic marine acquisition. For a perspective on dual sensor acquisition
devices and simultaneous measurements of pressure and vertical or normal
velocity, we refer to \cite{Carlson2007,Tenghamn2007}. Dual sensor
acquisition has additional benefits such as in noise reduction
\cite{Whitmore2010,Ronholt2015}. Independent of earlier analysis,
the results presented here confirm the relevance of acquiring dual
sensor data. 
The marine seismic acquisition we consider consists in sources positioned above 
the fixed receivers lattice. This setup matches the TopSeis 
acquisition system that has been recently deployed 
by CGG (Compagnie G\'en\'erale de G\'eophysique)
and Lundin Norway AS\footnote{see 
\url{https://www.cgg.com/en/What-We-Do/Offshore/Products-and-Solutions/TopSeis}
\label{footnote_topseis}},
for the purpose of improving the near offset data.

The formulation of FWI, and its underlying minimization procedure, 
has been proposed with norms different from energy norms. 
Shin and collaborators \cite{Shin2007a,Shin2007b,Shin2007c} 
compared the use of the phase and/or amplitude information in the data.
Envelope based misfit functional is considered in \cite{Bozdag2011,Wu2014}.
Brossier, Operto and Virieux investigate the use of the $L^1$ norm in 
\cite{Brossier2009}, and in combination with $L^2$ norm in \cite{Brossier2010}. 
In this paper, we introduce a new misfit functional which is based upon
conditional stability of the inverse problem for Cauchy data. This functional, 
related to the Green's identity, 
is formulated in terms of repeated integration of a quadratic 
expression. Such a formulation overcomes the difficulties of 
computational complexity occurring in discretizing operator norms, or 
distances of subspaces (which are typically used in theoretical 
stability estimates, \cite{Alessandrini2005,AdHGS}). 
Further, it allows independent locations for field and computational 
sources in the discrete settings. 

The paper is organized as follows. We detail the inverse 
problem and the main stability results in 
Section~\ref{section:stability_misfit}, where the misfit
functional to minimize for the reconstruction procedure 
is given. In Subsection~\ref{subsec:adjoint_state_misfit} 
the computation of the gradient of the misfit functional 
is conducted using the Lagrangian approach. In 
Section~\ref{section:numerical}, numerical experiments are
presented to demonstrate the efficiency of the algorithm, 
using single-frequency Cauchy data. Finally, 
Section~\ref{section:numerical_independent} shows
an experiment with different locations of observational 
and computational sources.

\section{Misfit functional and stability result}
\label{section:stability_misfit}
\subsection{Assumptions about the domain $\Omega$}\label{subsec assumption domain}

We fix some notations that will be adopted in this paper.
Given a point $x\in \mathbb{R}^3$, with $x=(x',x_3)$, where $x'\in\mathbb{R}^{2}$ and $x_3\in\mathbb{R}$, 
$B_r(x), B_r'(x')$ denote the open balls in
$\mathbb{R}^{3},\mathbb{R}^{2}$ centred at $x$, $x'$ respectively with radius $r$. We also denote by 
$Q_r(x)$ the cylinder

\[Q_r(x)=B_r'(x')\times(x_3-r,x_3+r)\]
and $B_r=B_r(0)$ and $Q_r=Q_r(0)$. 
\begin{enumerate}

\item We assume that $\Omega$ is a domain in $\mathbb{R}^3$
and that 
there exist positive constants $r_0$ and $B$ 
($r_0$ being dimensionally a length, $B$ being an absolute constant),
such that
\begin{equation}\label{assumption Omega}
|\Omega|\leq B r_0 ^3,
\end{equation}
where $|\Omega|$ denotes the Lebesgue measure of $\Omega$.


\item We fix an open non-empty subset $\Sigma$ of $\partial\Omega$
(where the measurements in terms of the local Cauchy data are taken).

\item We assume that $\Omega$ can be decomposed as follows:
\[\bar\Omega = \bigcup_{j=1}^{N}\bar{D}_j,\]
where $D_j$, $j=1,\dots , N$ are known open sets of
$\mathbb{R}^3$, satisfying the conditions below.

\begin{enumerate}
\item $D_j$, $j=1,\dots , N$ are connected and pairwise
non-overlapping polyhedrons.

\item $\partial{D}_j$, $j=1,\dots , N$ are of Lipschitz class with
constants $r_0$, $L$ see, for example, \cite{Adams2003}).

\item There exists one region, say $D_1$, such that
$\partial{D}_1\cap\Sigma$ contains a \emph{flat} portion
$\Sigma_1$ of size $r_0$ and for every 
$j\in\{2,\dots , N\}$,
we can select a subchain $\{ D_{j_k} \}_{k=1}^{K}$ ($K\leq N$)
of the partition $\{D_j\}_{j=1}^{N}$ of $\Omega$, such that
\begin{equation}\label{catena dominii}
D_{j_1}=D_1,\qquad D_{j_K}={D_j}.
\end{equation}

In addition, for every $k=1,\dots , K$,
${D}_{j_{k-1}}$ and ${D}_{j_{k}}$ are
contiguous in the sense that
$\partial{D}_{j_k}\cap \partial{D}_{j_{k-1}}$ contains a
\emph{flat} portion $\Sigma_k$ of size $r_0$ (here we agree that
$D_{j_0}=\mathbb{R}^3\setminus\Omega$), such that


\[\Sigma_k\subset\Omega,\quad\mbox{for\:every}\:k=2,\dots , K.\]
We illustrate the configuration and selection of subchain
in Figure~\ref{fig:sketch_geometry_1}.
We emphasize that under such an assumption, for every $k=1,\dots , K$, there exists $P_k\in\Sigma_k$ and
a rigid transformation of coordinates (depending on $k$) under which we have $P_k=0$
and

\begin{eqnarray}\label{flat}
\Sigma_k\cap{Q}_{r_{0}/3} &=&\{x\in
Q_{r_0/3}\:|\:x_3=0\},\nonumber\\
D_{j_k}\cap {Q}_{r_{0}/3} &=&\{x\in
Q_{r_0/3}\:|\:x_3>0\},\nonumber\\
D_{j_{k-1}}\cap {Q}_{r_{0}/3} &=&\{x\in
Q_{r_0/3}\:|\:x_3<0\}.
\end{eqnarray}

\end{enumerate}
\end{enumerate}

\newcommand{\geoblue} {blue!30!white}
\begin{figure}[ht!] \centering
\subfigure[Domain decomposition, the 
           ordering is arbitrary, at the 
           exception of $D_1$ which is 
           connected to $\Sigma$.]{
\begin{tikzpicture}[]
  \pgfmathsetmacro{\size}{4.25}
  \pgfmathsetmacro{\step}{5}  

  \draw[color=black,fill=yellow!20!white] (0,0)        rectangle (\size,\size);
  \draw[step=\size/\step.,gray] (0,0) grid (\size,\size);
  \draw[color=black] (0,0)        rectangle (\size,\size);

  \coordinate (corner) at (2*\size/\step,\step*\size/\step - \size/\step); 
  \draw (corner) node[anchor=south west,black,xshift=0cm]{$D_1$};  
  \coordinate (corner) at (3*\size/\step,\step*\size/\step - \size/\step); 
  \draw (corner) node[anchor=south west,black,xshift=0cm]{$D_2$};  
  \coordinate (corner) at (4*\size/\step,\step*\size/\step - \size/\step); 
  \draw (corner) node[anchor=south west,black,xshift=0cm]{$D_3$};  
  \coordinate (corner) at (0*\size/\step,\step*\size/\step - \size/\step); 
  \draw (corner) node[anchor=south west,black,xshift=0cm]{$D_4$};  
  \coordinate (corner) at (1*\size/\step,\step*\size/\step - \size/\step); 
  \draw (corner) node[anchor=south west,black,xshift=0cm]{$D_5$};  
  \coordinate (corner) at (0*\size/\step,\step*\size/\step - 2*\size/\step); 
  \draw (corner) node[anchor=south west,black,xshift=0cm]{$D_6$};  
  \coordinate (corner) at (1*\size/\step,\step*\size/\step - 2*\size/\step); 
  \draw (corner) node[anchor=south west,black,xshift=0cm]{$D_7$};  
  \coordinate (corner) at (2*\size/\step,\step*\size/\step - 2*\size/\step); 
  \draw (corner) node[anchor=south west,black,xshift=0cm]{$D_8$};  
  \coordinate (corner) at (3*\size/\step,\step*\size/\step - 2*\size/\step); 
  \draw (corner) node[anchor=south west,black,xshift=0cm]{$\dots$};  
    
  \coordinate (X0)at (1.5*\size/\step,\size);
  \coordinate (X1)at (3.5*\size/\step,\size);
  \draw[line width=3,red,]  (X0) to (X1);  

  \draw (X0) node[anchor=south west,red,xshift=0cm]{$\Sigma$};
  \draw (\size,0) node[anchor=south west,black,xshift=0cm]{$\Omega$};  

  \coordinate (X0)at (-0.05*\size,1.05*\size);
  \coordinate (X1)at ( 0.05*\size,1.05*\size);
  \coordinate (X2)at (-0.05*\size,0.95*\size);
  \draw[line width=1,gray,->] (X0) node[anchor=south west]{$x_1 x_2$} to (X1);
  \draw[line width=1,gray,->] (X0) node[anchor=east]      {$x_3$} to (X2);
    
\end{tikzpicture}
} \hspace*{1cm}
\subfigure[Path towards the selected subdomain, e.g. $D_{j_4}$. 
           $K=4$ in this illustration.]{
\begin{tikzpicture}[]
  \pgfmathsetmacro{\size}{4.25}  
  \pgfmathsetmacro{\step}{5}  
  
  \draw[color=black,fill=yellow!20!white] (0,0)        rectangle (\size,\size);
  \draw[step=\size/\step.,gray] (0,0) grid (\size,\size);
  \draw[color=black] (0,0)        rectangle (\size,\size);
  \draw[gray,fill=\geoblue] (3*\size/\step,2*\size/\step) rectangle (4*\size/\step,3*\size/\step);
  \draw (3*\size/\step,2*\size/\step) node[anchor=south west,black,xshift=0cm]
                                      {$D_{j_4}$};
  \draw (\size,0) node[anchor=south west,black,xshift=0cm]{$\Omega$};  
 
  \pgfmathsetmacro{\d}{\size/\step}  
  \coordinate(coo1) at (2*\d,\size-\d);
  \coordinate(coo2) at (3*\d,\size);
  \draw[gray,fill=\geoblue] (coo1) rectangle (coo2);
  \draw (coo1) node[anchor=south west,black,xshift=-1.5mm,yshift= 3mm]{$D_{1}=$ };
  \draw (coo1) node[anchor=south west,black,xshift= 2.0mm,yshift=-1mm]{$D_{j_1}$};
  \coordinate(coo1) at (2*\d,\size-2*\d);
  \coordinate(coo2) at (3*\d,\size-1*\d);
  \draw[gray,fill=\geoblue] (coo1) rectangle (coo2);
  \draw (coo1) node[anchor=south west,black,xshift=0cm]{$D_{j_2}$};  
  \coordinate(coo1) at (2*\d,\size-3*\d);
  \coordinate(coo2) at (3*\d,\size-2*\d);
  \draw[gray,fill=\geoblue] (coo1) rectangle (coo2);
  \draw (coo1) node[anchor=south west,black,xshift=0cm]{$D_{j_3}$};  

  \coordinate (X0)at (1.5*\size/\step,\size);
  \coordinate (X1)at (3.5*\size/\step,\size);
  \draw[line width=3,red,]  (X0) to (X1);  
  \draw (X0) node[anchor=south west,red,xshift=0cm]{$\Sigma$};
\end{tikzpicture}

}
\caption{Illustration of the domain 
         decomposition and the path towards
         a selected subdomain $D_{j_k}$.
         }
\label{fig:sketch_geometry_1}
\end{figure}
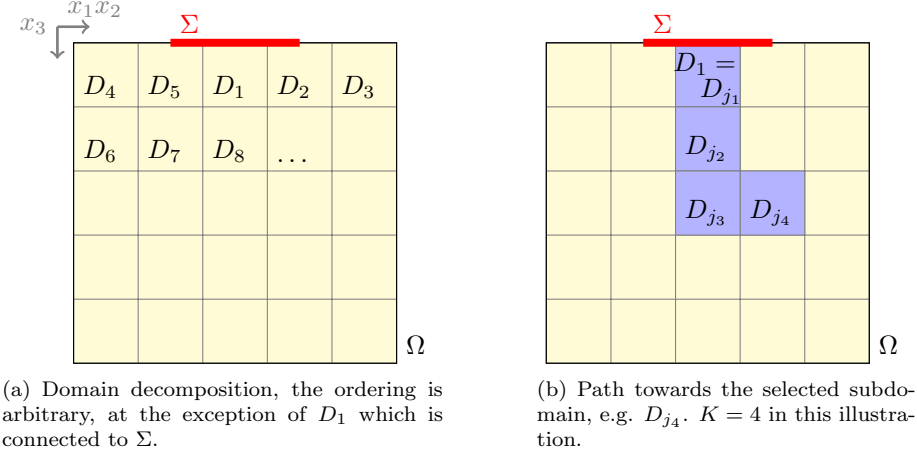


\subsection{A-priori information on the wave speed $c$}

We shall consider a real valued function $c\in L^{\infty}(\Omega)$, with

\begin{equation}\label{apriori q}
0<c_{min}\leq c\leq c_{max},
\end{equation}
for some positive constants $c_{min},\: c_{max}$ and of type

\begin{subequations}
\begin{eqnarray}\label{a priori info su q}
& &c(x)=\sum_{j=1}^{N}c_{j}(x)\chi_{D_j}(x),\qquad
x\in\Omega,\label{potential 1}\\
& & c_{j}(x)=a_j+A_j\cdot x\label{potential 2},
\end{eqnarray}
\end{subequations}
where $a_j\in\mathbb{R}$, $A_j\in\mathbb{R}^3$ 
are scalar and vector valued constants respectively,
and $\{D_j\}_{j=1}^{N}$, $j=1,\dots ,
N$ is the given partition of $\Omega$ introduced in 
Subsection~\ref{subsec assumption domain}.

\subsection{Cauchy data and misfit functional}
\label{subsection:cauchy_data_misfit}

We let $\Omega$ denote the subsurface domain with Lipschitz boundary
$\partial\Omega$ and $\Sigma$ the non-empty open portion of
$\partial\Omega$ introduced in Subsection~\ref{subsec assumption domain} where the acquisition is carried out. We introduce
the subspace of $H^{1/2}(\partial\Omega)$,
\[
   H^{1/2}_{co}(\Sigma) = \{f \in H^{1/2}(\partial\Omega) \:\vert\:
             \operatorname{supp} \, f \subset \Sigma\} .
\]
Its closure with respect to the $H^{1/2}(\partial\Omega)$ norm is the
space $H^{1/2}_{00}(\Sigma)$. In a similar manner, we define
$H^{1/2}_{00}(\partial\Omega \setminus \overline\Sigma)$. We denote
the pressure by $u$.
The \emph{Cauchy data} associated to $c$ is the space
$\mathcal{C}^{\Sigma}_c$,
\begin{eqnarray*}
   \mathcal{C}^{\Sigma}_c = \Big\{
   (f,g) \in H^{1/2}(\partial\Omega) |_{\Sigma} \times
             H^{-1/2}(\partial\Omega) |_{\Sigma}
         & \big| &\:
   \exists u \in H^{1}(\Omega)\ \text{weak solution to}\
   \Delta u + k^2 c^{-2} u = 0\ \text{in}\ \Omega,\\
  & & u \Big|_{\Sigma} =f,\\
   & & \langle \frac{\partial u}{\partial\nu} \Big|_{\partial\Omega},
       \varphi \rangle = \langle g, \varphi \rangle ,\
       \forall \varphi \in H^{1/2}_{00}(\Sigma) \Big\} .
\end{eqnarray*}
Here, $\langle\psi,\varphi\rangle$ denotes the duality between the complex 
valued spaces $H^{-\frac{1}{2}}(\partial\Omega)$, $H^{\frac{1}{2}}(\partial\Omega)$  
based on the $L^2$ inner product

\[\langle\psi,\varphi\rangle=\int_{\partial\Omega} \psi\overline\varphi\]
and $H^{\frac{1}{2}}(\partial\Omega)\big\vert_{\Sigma}$ and $H^{-\frac{1}{2}}(\partial\Omega)\big\vert_{\Sigma}$ 
denote the \emph{restrictions} of $H^{\frac{1}{2}}(\partial\Omega)$ and $H^{-\frac{1}{2}}(\partial\Omega)$ to 
$\Sigma$ respectively. $\mathcal{C}^{\Sigma}_c$ is a subspace of the Hilbert space 
$H^{\frac{1}{2}}(\partial\Omega)\big\vert_{\Sigma}\times H^{-\frac{1}{2}}(\partial\Omega)\big\vert_{\Sigma}$. 

We embed $\Omega$ in an ambient domain $\Upsilon \supset \Omega$ as 
we will find convenient to introduce Green's function not precisely 
for the physical domain $\Omega$ but for this augmented domain $\Upsilon$.

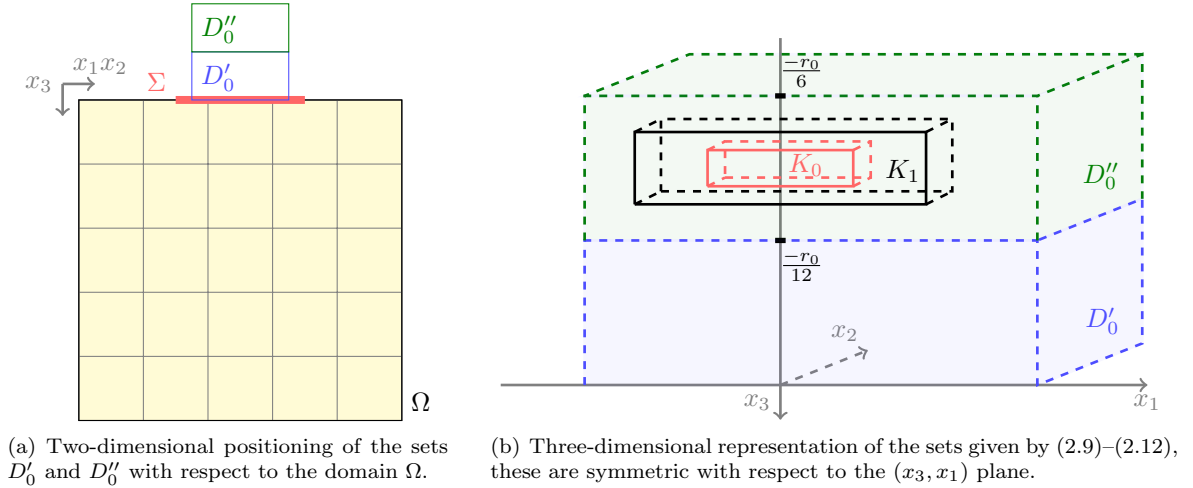
\begin{figure}[ht!]  \centering
\subfigure[Two-dimensional positioning of the 
           sets $D_0'$ and $D_0''$ with respect
           to the domain $\Omega$.]{
\begin{tikzpicture}[]
  \pgfmathsetmacro{\size}{4.25}
  \pgfmathsetmacro{\step}{5}  

  \draw[color=black,fill=yellow!20!white] (0,0)        rectangle (\size,\size);
  \draw[step=\size/\step.,gray] (0,0) grid (\size,\size);
  \draw[color=black] (0,0)        rectangle (\size,\size);
  \draw (\size,0) node[anchor=south west,black,xshift=0cm]{$\Omega$};  
 
  \coordinate (X0)at (1.5*\size/\step,\size);
  \coordinate (X1)at (3.5*\size/\step,\size);
  \draw[line width=3,\myred,]  (X0) to (X1);  
  \draw (X0) node[anchor=south east,\myred,xshift=0cm]{$\Sigma$};
  
  \coordinate (X0)at (1.75*\size/\step,\size);
  \coordinate (X1)at (3.25*\size/\step,1.15*\size);
  \draw[color=\myblue] (X0)        rectangle (X1);
  \draw (X0) node[anchor=south west,\myblue,xshift=0cm]{$D_0'$};

  \coordinate (X0)at (1.75*\size/\step,1.15*\size);
  \coordinate (X1)at (3.25*\size/\step,1.30*\size);
  \draw[color=\mygreen] (X0)        rectangle (X1);    
  \draw (X0) node[anchor=south west,\mygreen,xshift=0cm]{$D_0''$};

  \coordinate (X0)at (-0.05*\size,1.05*\size);
  \coordinate (X1)at ( 0.05*\size,1.05*\size);
  \coordinate (X2)at (-0.05*\size,0.95*\size);
  \draw[line width=1,gray,->] (X0) node[anchor=south west]{$x_1 x_2$} to (X1);
  \draw[line width=1,gray,->] (X0) node[anchor=east]      {$x_3$} to (X2);
  
\end{tikzpicture}
} \hfill
\subfigure[Three-dimensional representation 
           of the sets given by 
           \eqref{eq:D0prime}--\eqref{eq:K1},
           these are symmetric with respect to the 
           $(x_3,x_1)$ plane.]{
\begin{tikzpicture}[]
  \pgfmathsetmacro{\size}{2.30}
  \pgfmathsetmacro{\ratiominx}{-1.60}
  \pgfmathsetmacro{\ratiomaxx}{ 2.10}
  \pgfmathsetmacro{\ratiominz}{-0.20}  
  \pgfmathsetmacro{\ratiomaxz}{ 2.00}
  \pgfmathsetmacro{\ratioyx}  { 0.50}  
  \pgfmathsetmacro{\ratioyz}  { 0.20}
  \coordinate (O)   at (         0, 0); 
  \coordinate (Xmin)at (\ratiominx*\size,0);
  \coordinate (Xmax)at (\ratiomaxx*\size,0);
  \coordinate (Zmin)at (         0      ,\ratiominz*\size);
  \coordinate (Zmax)at (         0      ,\ratiomaxz*\size);
  \coordinate (Ymax)at (\ratioyx*\size  ,\ratioyz*\size);
  \draw[line width=1,gray,->]  (Xmin) to (Xmax);  
  \draw[line width=1,gray,<-]  (Zmin) to (Zmax);  
  \draw[line width=1,gray,->,dashed]     (O) to (Ymax);  
  \draw (Zmin) node[anchor=south east,gray,xshift=0cm]{$x_3$};
  \draw (Xmax) node[anchor=north east,gray,xshift=3mm]{$x_1$};
  \draw (Ymax) node[anchor=south east,gray,xshift=0cm]{$x_2$};
  \pgfmathsetmacro{\scale}{10}
  \pgfmathsetmacro{\lim}  {0.70}
  \pgfmathsetmacro{\dmin}{0} 
  \pgfmathsetmacro{\dmax}{\scale*1./12.}
  \coordinate (D11) at (\lim*\ratiominx*\size,\dmin*\size);
  \coordinate (D12) at (\lim*\ratiomaxx*\size,\dmin*\size);  
  \coordinate (D13) at (\lim*\ratiomaxx*\size,\dmax*\size);  
  \coordinate (D14) at (\lim*\ratiominx*\size,\dmax*\size); 
  \pgfmathsetmacro{\shiftx}  {1.20*\ratioyx*\size}
  \pgfmathsetmacro{\shiftz}  {1.20*\ratioyz*\size}
  \coordinate (D21) at (\shiftx + \lim*\ratiominx*\size,\shiftz + \dmin*\size);
  \coordinate (D22) at (\shiftx + \lim*\ratiomaxx*\size,\shiftz + \dmin*\size);  
  \coordinate (D23) at (\shiftx + \lim*\ratiomaxx*\size,\shiftz + \dmax*\size);  
  \coordinate (D24) at (\shiftx + \lim*\ratiominx*\size,\shiftz + \dmax*\size);
  \draw[line width=1,\myblue,dashed] 
                     (D11) to (D14) to (D13) to (D12) to (D22) to
                     (D23) to (D13)    ;
  \draw[line width=1,\myblue,fill=\myblue,opacity=.05] 
            (D12) to (D11) to (D14) to (D13) to (D12) to (D22) to
                     (D23) to (D13)    ;
  \draw (D22) node[anchor=south east,\myblue,xshift=-.5em]{$D_0'$};
  \coordinate (x-) at ( -.03*\size,\dmax*\size);
  \coordinate (x+) at (  .03*\size,\dmax*\size);
  \draw[line width=2,black] (x-) to (x+);
  \draw (x+) node[anchor=north west,black,xshift=-.2cm]{$\frac{-r_0}{12}$};
    
  \pgfmathsetmacro{\scale}{10}
  \pgfmathsetmacro{\lim}  {0.70}
  \pgfmathsetmacro{\dmin}{\scale*1./12.} 
  \pgfmathsetmacro{\dmax}{\scale*1./6.}
  \coordinate (D11) at (\lim*\ratiominx*\size,\dmin*\size);
  \coordinate (D12) at (\lim*\ratiomaxx*\size,\dmin*\size);  
  \coordinate (D13) at (\lim*\ratiomaxx*\size,\dmax*\size);  
  \coordinate (D14) at (\lim*\ratiominx*\size,\dmax*\size); 
  \pgfmathsetmacro{\shiftx}  {1.20*\ratioyx*\size}
  \pgfmathsetmacro{\shiftz}  {1.20*\ratioyz*\size}
  \coordinate (D21) at (\shiftx + \lim*\ratiominx*\size,\shiftz + \dmin*\size);
  \coordinate (D22) at (\shiftx + \lim*\ratiomaxx*\size,\shiftz + \dmin*\size);  
  \coordinate (D23) at (\shiftx + \lim*\ratiomaxx*\size,\shiftz + \dmax*\size);  
  \coordinate (D24) at (\shiftx + \lim*\ratiominx*\size,\shiftz + \dmax*\size);
  \draw[line width=1,\mygreen,dashed] 
                     (D11) to (D14) to (D13) to 
                     (D23) to (D22) ;
  \draw[line width=1,\mygreen,dashed] 
                     (D23) to (D24) to (D14) ;
  \draw[line width=1,\mygreen,dashed] 
                     (D12) to (D13) ;
  \draw[line width=1,\myblue,fill=\mygreen,opacity=.05] 
                     (D23) to (D24) to (D14) to (D13) ;
  \draw[line width=1,\myblue,fill=\mygreen,opacity=.05] 
            (D12) to (D11) to (D14) to (D13) to (D12) to (D22) to
                     (D23) to (D13)    ;
  \draw (D22) node[anchor=south east,\mygreen,xshift=-.5em]{$D_0''$};
  \coordinate (x-) at ( -.03*\size,\dmax*\size);
  \coordinate (x+) at (  .03*\size,\dmax*\size);
  \draw[line width=2,black] (x-) to (x+);
  \draw (x+) node[anchor=south west,black,xshift=-.2cm]{$\frac{-r_0}{6}$};
  
  \pgfmathsetmacro{\dzmin}{\scale*11./96.}
  \pgfmathsetmacro{\dzmax}{\scale*13./96.}
  \pgfmathsetmacro{\dxmin}{-\scale*1./24.}
  \pgfmathsetmacro{\dxmax}{ \scale*1./24.}
  \pgfmathsetmacro{\dymin}{-\scale*1./24.}
  \pgfmathsetmacro{\dymax}{ \scale*1./24.}
  \coordinate (D11) at (\dxmin*\size,         \dzmin*\size);
  \coordinate (D12) at (\dxmax*\size,         \dzmin*\size);  
  \coordinate (D13) at (\dxmax*\size,         \dzmax*\size);  
  \coordinate (D14) at (\dxmin*\size,         \dzmax*\size);  

  \coordinate (D21) at (\dxmin*\size+0.1*\size,.05*\size+\dzmin*\size);
  \coordinate (D22) at (\dxmax*\size+0.1*\size,.05*\size+\dzmin*\size);  
  \coordinate (D23) at (\dxmax*\size+0.1*\size,.05*\size+\dzmax*\size);  
  \coordinate (D24) at (\dxmin*\size+0.1*\size,.05*\size+\dzmax*\size);
  \draw[line width=1,\myred] 
              (D11) to (D12) to (D13) to (D14) to (D11);
  \draw[line width=1,\myred,dashed] 
              (D21) to (D22) to (D23) to (D24) to (D21);
  \draw[line width=1,\myred,dashed] (D12) to (D22);
  \draw[line width=1,\myred,dashed] (D13) to (D23);
  \draw[line width=1,\myred,dashed] (D14) to (D24);
  \draw[line width=1,\myred,dashed] (D11) to (D21);
  \draw (D22) node[anchor=south east,\myred,xshift=-.5cm,yshift=-1mm]{$K_0$};

  \pgfmathsetmacro{\dzmin}{ \scale*5./48.}
  \pgfmathsetmacro{\dzmax}{ \scale*7./48.}
  \pgfmathsetmacro{\dxmin}{-\scale*1./12.}
  \pgfmathsetmacro{\dxmax}{ \scale*1./12.}
  \pgfmathsetmacro{\dymin}{-\scale*1./12.}
  \pgfmathsetmacro{\dymax}{ \scale*1./12.}
  \coordinate (D11) at (\dxmin*\size,         \dzmin*\size);
  \coordinate (D12) at (\dxmax*\size,         \dzmin*\size);  
  \coordinate (D13) at (\dxmax*\size,         \dzmax*\size);  
  \coordinate (D14) at (\dxmin*\size,         \dzmax*\size);  

  \coordinate (D21) at (\dxmin*\size+0.15*\size,.075*\size+\dzmin*\size);
  \coordinate (D22) at (\dxmax*\size+0.15*\size,.075*\size+\dzmin*\size);  
  \coordinate (D23) at (\dxmax*\size+0.15*\size,.075*\size+\dzmax*\size);  
  \coordinate (D24) at (\dxmin*\size+0.15*\size,.075*\size+\dzmax*\size);
  \draw[line width=1,black] 
              (D11) to (D12) to (D13) to (D14) to (D11);
  \draw[line width=1,black,dashed] 
              (D21) to (D22) to (D23) to (D24) to (D21);
  \draw[line width=1,black,dashed] (D12) to (D22);
  \draw[line width=1,black,dashed] (D13) to (D23);
  \draw[line width=1,black,dashed] (D14) to (D24);
  \draw[line width=1,black,dashed] (D11) to (D21);
  \draw (D22) node[anchor=south east,black,xshift=-3mm]{$K_1$};
    
\end{tikzpicture}
}
\caption{Illustration of the sets 
        that are defined in 
        \eqref{eq:D0prime}--\eqref{eq:K1}.}
\label{fig:sketch_geometry_2}
\end{figure}

We recall that by assumption 3(c) of Subsection~\ref{subsec assumption domain} we can assume 
that there exists a point $P_1$ such that up to a rigid transformation of coordinates we have 
that $P_1=0$ and \eqref{flat} holds with $\Sigma=\Sigma_1$. Denoting by

\[D_0=\left\{x\in(\mathbb{R}^3\setminus\Omega)\cap B_{\frac{r_0}{3}}\:\bigg|\:|x_i|<\frac{r_0}{6},\:i=1,2\:; \:-\frac{r_0}{6}< x_3<0\right\},\]
it turns out that the augmented domain $\Upsilon=\accentset{\circ}{\overline{(\Omega\cup D_0)}}$ is of Lipschitz class with constants $\frac{r_0}{3}$ and $\widetilde{L}$, where $\widetilde{L}$ depends on $L$ only.
Given $r>0$, we set

\begin{eqnarray}
&&\Gamma_1=\left\{ x\in \Upsilon \ \bigg|\: |x_i|< \frac{r_0}{6},\ i=1,\ 2\  ; \  x_3=-\frac{r_0}{6} \right\}, \\
&&\Gamma_2=\partial\Upsilon\setminus\overline\Gamma_1,\\
&&(\Upsilon)_r = \left\{x\in \Upsilon\ \bigg| \ \mbox{dist}(x,\partial \Upsilon)>r \right\}.
\end{eqnarray}
We also introduce the following sets

\begin{eqnarray}
  &&  D^{'}_0=\left\{x\in D_0\:\bigg|\:-\frac{r_0}{12}< x_3 <0\right\}, 
     \label{eq:D0prime}    \\
  &&  D^{''}_0=\left\{x\in D_0\:\bigg|\:-\frac{r_0}{6}< x_3 <-\frac{r_0}{12}\right\}, 
     \label{eq:D0second}\\
  &&  K_0=\left\{x\in D^{''}_0\:\bigg|\:|x_i|<\frac{r_0}{24},
     \:i=1,2\:;\:-\frac{13}{96}r_0< x_3 <-\frac{11}{96}r_0\right\}, 
     \label{eq:K0}\\
  &&  K_1=\left\{x\in D^{''}_0\:\bigg|\:|x_i|<\frac{r_0}{12},
      \:i=1,2\:;\:-\frac{7}{48}r_0< x_3 <-\frac{5}{48}r_0\right\},
      \label{eq:K1} \\
  && \Omega'=\accentset{\circ}{\overline{(D^{'}_0\cup\Omega)}}
      \label{eq:Omegaprime}.
\end{eqnarray}
Note that, fixing the origin at the center of $D^{''}_0, K_0,K_1$ are 
concentric parallelograms 
scaled by the  factors $\frac{1}{4}, \frac{1}{2}$ respectively.
We illustrate the geometry in Figure~\ref{fig:sketch_geometry_2}.
Note also that this precise choice of scale parameters is just made for the sake of 
definiteness. What really matters is the general geometrical configuration, in particular 
we must have

\[\mbox{dist}(D^{''}_0,\overline{\Omega})\geq \frac{r_0}{12}>0\quad\textnormal{and}\quad K_0\subset\subset K_1\subset\subset D^{''}_0.\]
We shall denote by $\Gamma(x,y)$ the standard fundamental solution to the Laplace equation which is
\begin{eqnarray}
\Gamma(x,y)=\frac{1}{4\pi |x-y|} \ .
\end{eqnarray}

\begin{definition}
Let $B$, $N$, $r_0$, $L$, $c_{min}$ and $c_{max}$ be given
positive numbers with $N\in\mathbb{N}$. We will
refer to this set of numbers as to the \textit{a-priori data}. Several constants depending on 
the \emph{a-priori data} will appear within the paper. In order 
to simplify our notation, any quantity denoted by $C,C_1,C_2, \dots$ 
will be called a \emph{constant} understanding in most cases that 
it only depends on the \textit{a-priori} data. 
\end{definition}

Next we introduce a mixed boundary value problem for $\Delta + k^2c^{-2}$ in 
$\Upsilon$ which is always well posed, independently of any \emph{a-priori} 
condition on $c$, besides the assumption of being real valued and bounded. 
This shall enable us to construct Green's function for $\Delta + k^2c^{-2}$ 
in $\Upsilon$.

We assume that the uppermost part $D_0$ of the domain 
$\Upsilon$ represents a region filled by water. The 
wave speed in $D_0$ can then be assumed to be known
and constant with 
\begin{equation}
  c = c_0, \quad\text{in $D_0$.}
\end{equation}
Pressure sources (air gun) are excited to produce 
impulses located at points in $K_1 \subset D_0$ and 
Cauchy data are collected through dual sensors located 
at the surface $\Sigma$, which lies below $K_1$.

Pressure is assumed to be zero at the sea level $\Gamma_1$
(\textit{i.e.}, free surface), and to satisfy, on the remaining
part of the boundary $\Gamma_2$ (of the region of interest), 
a (conventional) absorbing condition.
If we model this problem in the frequency domain, and 
assume that the source is modeled by a Dirac's delta 
concentrated at a point $y$, the pressure is represented 
as the Green's function of the following mixed boundary 
value problem,
\begin{equation}\label{Green_third_kind_statement}
\left\{
\begin{array}
{lcl}  \Delta \Gc (\cdot, y) +k^2c^{-2}(\cdot) \Gc(\cdot, y) =-\delta(\cdot-y)\ ,&&
\mbox{in $\Upsilon$ ,}
\\
\Gc(\cdot,y)= 0\ ,   && \mbox{on $\Gamma_1$ ,}
\\
\partial _{\nu}{\Gc}(\cdot,y) - i k_0 \Gc(\cdot, y) =0 \ , && \mbox{on $\Gamma_2$ .}
\end{array}
\right.
\end{equation}
The theory developed, for example, in 
\cite{Cakoni2014} shows that such a function 
$\Gc$ exists and is unique 
in the case of constant wave speed $c$. 
Note that the term $k_0$ is conventionally assumed to be 
constant and known.
The next theorem collects the main 
features of the Green's function solving
\eqref{Green_third_kind_statement} 
also in the case of variable wave speed $c(\bx)$.
A similar result, but with stronger hypothesis, was proven
in \cite[Propositions~3.1,~3.4,~3.5]{AdHGS}. The thesis 
here is slightly weaker, but the argument is somewhat 
simpler.

\begin{theorem}\label{wellposedness} 
For any $y \in \Upsilon$, there exists a unique 
distributional solution $\Gc(\cdot,y)$ to 
\eqref{Green_third_kind_statement}. Moreover, there
exists a constant
$C>0$ depending on $r_0,L, k$  and on $c_{min}$ such that 
for any $x,y\in {(\Upsilon)}_{r_0}$, $x\neq y$ we have that 
\begin{eqnarray}\label{1}
|\Gc(x,y) - \Gamma(x,y)|\le C \ 
\end{eqnarray}
and
\begin{eqnarray}\label{2}
|{\nabla}_y \Gc(x,y) - {\nabla}_y \Gamma(x,y)|\le C (\big| \log|x-y| \big| +1) 	\ .
\end{eqnarray}
Moreover, let $Q_{k+1}$ be a point such that $Q_{k+1}\in B_{\frac{r_0}{8}}(P_{k+1})\cap 
\Sigma_{k+1}$ with $k\in \{1, \dots, N-1 \}$, then the following inequality holds true for 
every $x\in B_{\frac{r_0}{16}}(P_{k+1})\cap D_{j_k+1}$ and every $y=Q_{k+1}-re_3$, where 
$r\in (0,\frac{r_0}{16})$

\begin{eqnarray}\label{3}
|{\nabla}^2_y(\Gc(x,y)-\Gamma(x,y))|\le Cr^{-1}.
\end{eqnarray}
\end{theorem}
Here ${\nabla}^2_y$ denotes the Hessian matrix.

\begin{proof}
Let $y\in \Upsilon$ and let $G_0$ be the Green's function for the Laplace 
operator which solves 
\begin{equation}\label{Green_Laplace_1}
\left\{
\begin{array}
{lcl}  \Delta G_0 (\cdot, y) =-\delta(\cdot-y)\ ,&&
\mbox{in $\Upsilon$ ,}
\\
 G_0(\cdot,y)= 0\ ,   && \mbox{on $\Gamma_1$ ,}
\\
\partial _{\nu}{G_0}(\cdot,y) - i k_0 G_0(\cdot, y) =0 \ , && \mbox{on $\Gamma_2$ .}
\end{array}
\right.
\end{equation}
The existence and uniqueness of a distributional solution $G_0\in L^1(\Upsilon)$ 
to \eqref{Green_Laplace_1} is a consequence of standard theory on boundary value 
problems for the Laplace equation. 
By standard techniques it can be proved that for any $y\in \Upsilon$ such that 
$\mbox{dist}(y,\partial \Upsilon)\ge \frac{r_0}{4}$ we have that 
$G_0(\cdot,y)\in L^2(\Upsilon)$. Now we define $\Rc(\cdot,y)\in H^1(\Upsilon)$  to be
the solution to 

\begin{equation}\label{Schrodinger_rhs}
\left\{
\begin{array}
{lcl}  \Delta \Rc (\cdot, y) + k^2 c^{-2}\Rc(\cdot,y)= -k^2 c^{-2} G_0(\cdot,y)\ ,&&
\mbox{in $\Upsilon$ ,}
\\
 \Rc(\cdot,y)= 0\ ,   && \mbox{on $\Gamma_1$ ,}
\\
\partial _{\nu}{\Rc}(\cdot,y) - i k_0  \Rc(\cdot, y) =0 \ , && \mbox{on $\Gamma_2$ .}
\end{array}
\right.
\end{equation}
The existence and uniqueness for \eqref{Schrodinger_rhs} follows along the lines of 
the proof of \cite[Proposition 3.1]{AdHGS}, which relies on the Fredholm altenative 
theory.
Moreover, by arguments based on well-known estimates for the Cauchy problem contained 
in \cite[Proposition 3.1]{AdHGS}, we have that
\begin{eqnarray}
\|\Rc(\cdot,y) \|_{H^1(\Upsilon)}\le C \|G_0(\cdot,y) \|_{L^2(\Upsilon)}
\end{eqnarray}
and, by standard interior estimates, that

\begin{eqnarray}
|\Rc(x,y)|\le C ,
\end{eqnarray}
for any $x\neq y, \ x\in \Upsilon $ and $\mbox{dist}(y,\partial\Upsilon)\ge \frac{r_0}{4}$.
If we form 

\begin{equation} \label{eq:GR}
\Gc(x,y)=G_0(x,y) +\Rc(x,y),
\end{equation}
then we end up with the following estimate 

\begin{eqnarray}
|\Gc(x,y)|\le C |x-y|^{-1},
\end{eqnarray}
for any $x,y\in \Upsilon,\ x\neq y$ and $\mbox{dist}(x,\partial\Upsilon)\ge \frac{r_0}{4}, 
\:\mbox{dist}(y,\partial\Upsilon)\ge \frac{r_0}{4}$ . The latter combined with the arguments 
in the proofs of \cite[Propositions 3.4, 3.5]{AdHGS} suffices to deduce \eqref{1}, \eqref{2} and \eqref{3}. 
\end{proof}

Assuming that, for sources placed at arbitrary 
points $z \in K_1$, we can measure associated 
Cauchy data on $\Sigma$:
\begin{equation}
  \Gobs(x,z),~\dfrac{\partial}{\partial \nu_x}\Gobs(x,z),\quad x\in\Sigma,
\end{equation}
we seek $c$ which minimizes the 
following misfit functional

\begin{equation}\label{misfit}
\mathcal{J}(c) = \int_{K_1 \times K_1} 
                 \left|\int_{\Sigma} 
                 \left(\Gc(x,y)\:\partial_{\nu}\Gobs(x,z)
                      -\Gobs(x,z)\:\partial_{\nu}\Gc(x,y)\right)
                      \dd \mu(x) \right|^2\:\dd y\:\dd z,
\end{equation}
where $\mu$ denotes the element of surface measure. \\

The introduction of the misfit functional~\eqref{misfit} is motivated 
by the following argument. Given two wave speeds 
$c^{(i)}$, $i=1,2$, consider the Green's functions $G_i$, introduced in 
Theorem~\ref{wellposedness}, corresponding to $c^{(i)}$ in $\Upsilon$ and the following quantity

\begin{equation}\label{J}
\mathcal{J}(c^{(1)},c^{(2)}) =\int_{K_1 \times K_1} 
                              \left|\mathcal{S}_{\mathcal{U}_0}(y,z)\right|^2 \dd y\:\dd z,
\end{equation}
where
\begin{equation}\label{SU 0}
\mathcal{S}_{\mathcal{U}_0}(y,z)=\int_{\Sigma_1} \left(G_1(x,y)\:\partial_{\nu}G_2(x,z)-G_2(x,z)\:\partial_{\nu}G_1(x,y)\right) \dd \mu(x),\quad\mbox{for\:any}\:y,z\in D_0.
\end{equation}
Expressions of the form above have appeared in many occasions
in the treatment of inverse boundary problems. Analogies can be 
found with the probe method by Ikehata \cite{Ikehata1998}, see
also \cite{Potthast2006}. In particular, a very strong relation
can be observed with the so-called \emph{reciprocity gap functional}
introduced by Colton and Haddar \cite{ColtonHaddar2005} for 
inverse scattering.

It would be a matter of an exercise to show that,
fixing $y\in K_1$,
\begin{equation}
  c \in L^\infty(\Upsilon) \rightarrow 
                           \big(G_c(x,y)|_\Sigma,~\partial_\nu G_c(x,y)|_\Sigma \big)
                           \in H^{1/2}(\partial\Omega)|_\Sigma \times
                               H^{-1/2}(\partial\Omega)|_\Sigma
\end{equation}
is Fr\'echet differentiable. Note also that, since we are 
assuming~\eqref{potential 1} and~\eqref{potential 2} 
(that is, $c$ lives in a finite dimensional space), the 
$L^\infty(\Upsilon)$ norm can equivalently be replaced by the 
$L^2(\Upsilon)$ norm. This will enable us to apply to $\misfit$ a 
projected steepest descent method in Section~\ref{section:numerical}.

\begin{theorem}\label{teorema principale}
Let $\Omega$, $D_j$, $j=1,\dots , N$ and $\Sigma$ be a domain, $N$ subdomains of $\Omega$ and a portion of $\partial\Omega$ as in section \ref{subsec assumption domain} respectively.
Let $c^{(i)}$, $i=1,2$ be two wave speeds satisfying \eqref{apriori q} and of type

\begin{equation}\label{a priori info su sigma}
c^{(i)}=\sum_{j=1}^{N}c^{(i)}_{j}(x)\chi_{D_j}(x),\qquad
x\in\Omega,
\end{equation}
where
\[c^{(i)}_{j}(x)=a^{(i)}_j+A^{(i)}_j\cdot x,\]
with $a^{(i)}_j\in\mathbb{R}$ and $A^{(i)}_j\in\mathbb{R}^3$, then we have
\begin{equation}\label{stabilita' globale}
||c^{(1)}-c^{(2)}||_{L^{\infty}(\Omega)}\leq C \left(\mathcal{J}(c^{(1)},c^{(2)})\right)^{1/2},
\end{equation}
where $C$ is a positive constant that depends on the a-priori data
only.
\end{theorem}

\begin{remark}
Note that the introduction of the misfit functional is 
driven on the one hand by our computational experiments and 
on the other hand, it is inspired by the method of 
singular solutions used in previous stability results (see \cite{AdHGS} for the case of the Helmholtz equation). Although 
a natural metric on the space of Cauchy data 
$\mathcal{C}_{c}^{\Sigma}$ is given by the 
\emph{distance (aperture)} $d$ introduced in \cite[(2.8)]{AdHGS}, 
we have

\begin{equation}\label{J and C}
\mathcal{J}(c^{(1)},c^{(2)})\leq C d(\mathcal{C}^{\Sigma}_{c^{(1)}}, \mathcal{C}^{\Sigma}_{c^{(2)}})^2
\end{equation}
and also 
\begin{equation}\label{J and Impedance}
\mathcal{J}(c^{(1)},c^{(2)})\leq C ||\Lambda^\Sigma_{c^{(1)}} - \Lambda^\Sigma_{c^{(2)}}||_{*}^{2},
\end{equation}
where  
$\Lambda^\Sigma_c: u \in H^{1/2}(\partial\Omega)|_\Sigma \rightarrow
 \partial_\nu u \in H^{-1/2}(\partial\Omega)|_\Sigma $ 
 is the local Dirichlet-to-Neumann map
with its natural norm (here denoted by $||\cdot||_{*}$) between local trace spaces. These estimates justify the use of $\mathcal{J}$ as a substitute to the more traditional quantifications of the error on boundary data (either involving distance of spaces of Cauchy data or boundary maps).
\end{remark}

\begin{proof}{\textit{of Theorem \ref{teorema principale}}}.

The proof requires only some adaptations of 
\cite[Theorem 2.2]{AdHGS}, which are outlined below.

\begin{description}
\item [i)] We introduce different boundary conditions. This aspect involves some modifications in the constructions of the Green's function and has been treated in Theorem $1$ above. Note that we took advantage of the fact that here we focus on the three-dimensional case only.

\item [ii)] We replace the domain $\Omega$ with 
$\Omega '$ and add as initial subdomain $D_0'$ 
instead of $D_1$. We take advantage of the fact 
that $c^{(1)}=c^{(2)}=1$ in 
$D_0'\subset \Upsilon\setminus \overline{\Omega}$ 
which allows us to skip the arguments 
of \cite[Section 4.3]{AdHGS}. 

\item [iii)] We observe that in \cite[Theorem 2.2]{AdHGS} 
             the right-hand side in formula $(2.20)$ is 
             expressed in terms of the distance between spaces
             of Cauchy data. However, the only Cauchy data 
             that are actually used are those arising from
             Green's function with pole
             in $\widehat{D}_0$.
             And the role of $\widehat{D}_0$ 
             (see \cite{AdHGS}) can be 
             equivalently taken by the set $K_0$ introduced 
             here in \eqref{eq:K0}. Moreover, such Cauchy
             data intervene only in expressions like the one
             in~\eqref{SU 0} above. Therefore, the right-hand
             side of~\cite[(2.20)]{AdHGS} can be replaced by
\begin{eqnarray}
\sup _{K_0\times K_0}|S_{\mathcal{U}_0}(y,z)|.
\end{eqnarray}

We also recall that $S_{\mathcal{U}_0}(y,z)$ is a solution to 
\begin{eqnarray}
(\Delta_y +\Delta_z +2k^2)S_{\mathcal{U}_0}(y,z)=0
\end{eqnarray} 
in $D_0\times D_0$ (see \cite[(4.61)]{AdHGS}). Consequently, by standard estimate of boundedness in the interior we have 
\begin{eqnarray}
\sup_{K_0\times K_0} |S_{\mathcal{U}_0}(y,z)|^2\le C \int_{K_1\times K_1} |S_{\mathcal{U}_0}(y,z)|^2\:\dd y\:\dd z
\end{eqnarray}
where $C>0$ is a constant depending on $k$ and on $r_0$ . 
\item  [iv)]

Another difference comes from the fact that we are now assuming $c^{(i)}$ piecewise linear instead of $q^{(i)}=k^2 (c^{(i)})^{-2}$, for $i=1,2$. However due to the assumption \eqref{apriori q} the estimation of $c^{(i)}, \nabla c^{(i)}$ at each interface is equivalent to that for $q^{(i)}, \nabla q^{(i)}$, for $i=1,2$.
\end{description}

\end{proof}

The stability given in Theorem~\ref{teorema principale} 
justifies the use of \eqref{misfit} for an optimization
algorithm toward the inversion of \emph{Cauchy data to 
wave speed}.

\subsection{Computation of the gradient using the adjoint state method}
\label{subsec:adjoint_state_misfit}
We start by observing that, although $c\rightarrow \Gc(x,y)$ 
does not map into $H^1(\Upsilon)$ (because of the singularity
of the Green's function), the derivative $\dGc(x,y)[c]\dc$ 
exists and does belong in $H^1(\Upsilon)$. This can be 
achieved by recalling~\eqref{Schrodinger_rhs} and~\eqref{eq:GR}. That 
is 
\begin{equation}
   \Rc(x,y) = \Gc(x,y) - G_0(x,y) \quad \in H^1(\Upsilon),
\end{equation}
and the second term $G_0$ is independent of $c$.
Hence we may well define
\begin{equation}
  \dGc(x,y)[c]\dc =  D_c\Rc(x,y)[c]\dc.
\end{equation}
We denote
\begin{equation}
  H_{\Gamma_1}^1(\Upsilon) = \{v \in H^1(\Upsilon) ~|~ v|_{\Gamma_1}=0 \}.
\end{equation}

\newcommand{\vbar}{v}
We continue with the variational formulation of 
Problem~\eqref{Green_third_kind_statement}.
As noted above, denoting $\Rc(x,y)=G_c(x,y)-G_0(x,y)$,
it can be formulated as 
\begin{equation} \label{variational_formulation}
    \int_\Upsilon k^2 c^{-2}(x) \Gc(x,y) \vbar(x) 
  - \nabla_x \Rc(x,y) \cdot \nabla_x \vbar(x)\:\dd x 
  + \int_{\Gamma_2} \ii k_0 \Rc(x,y) \vbar(x)\:\dd \mu(x) = 0,
\end{equation}
for every $v \in H_{\Gamma_1}^1(\Upsilon)$.

The parameter reconstruction is conducted via an iterative 
minimization of the misfit functional $\misfit$ of~\eqref{misfit}, 
in a gradient descent algorithm. Therefore, we require the 
computation of the gradient of $\misfit$. 
For this purpose, we employ the adjoint state method, which 
allows the computation of the gradient without having to 
form explicitly the derivative $\Gc'$. The method arose from 
the work of \cite{Lions1971} and was promoted in the context 
of parameter derivation in \cite{Chavent1974}. It has massively 
been employed since then, and we refer to \cite{Plessix2006}
for a review in the geophysical framework.
Here, we follow the traditional steps 
for the selection of the Lagrange multiplier and 
formation of the gradient which are detailed, 
for example, in \cite{Chavent2010,Kern2016}, 
and that we adapt to our choice of misfit functional.

We first postpone the sum over the sources in the 
misfit functional~\eqref{misfit}, and select a single source 
for $\Gc$ and $\Gobs$, $y$ and $z$ respectively. We introduce
\begin{equation} \label{misfit_unisource}
  \misfitI(c)(y,z) = \left|\int_{\Sigma} \Big(\Gc(x,y)\:\partial_{\nu}
                    \Gobs(x,z)-\Gobs(x,z)\:\partial_{\nu}\Gc(x,y)\Big)
                    \dd \mu(x) \right|^2
                   = \left| \Su(y,z) \right|^2,
\end{equation}
such that
\begin{equation} \label{equation:misfit_I}
  \misfit(c) = \int_{K_1 \times K_1} \misfitI(c)(y,z) \:\dd y\:\dd z.
\end{equation}
The Riesz representation theorem gives 
\begin{equation}
  D_c \misfitI[c] \dc = \int_\Upsilon \nabla_c \misfitI \dc\:\dd x,
\end{equation}
where $\nabla_c\misfitI$ is the gradient, and 
$D_c \misfitI$ is the differential defined 
for every $\dc \in L^2(\Upsilon)$ by
\begin{equation}
  D_c \misfitI[c] ~:~ \dc \quad \rightarrow \quad 
                 \lim_{h\rightarrow 0} \dfrac{\misfitI(c + h\dc) - \misfitI(c)}{h}.
\end{equation}

The adjoint state method considers 
the constrained minimization problem
\begin{equation}
\begin{aligned}
 \min_c \misfitI(c) \quad \quad 
 \text{subject to~\eqref{Green_third_kind_statement}}.
\end{aligned}
\end{equation}
The constraint can be replaced by the variational
formulation \eqref{variational_formulation} and
the associated formulation of the Lagrangian is defined by
\begin{equation}
\begin{aligned}
  \lagrangian(c,G,\tildadj) = 
     \misfitI(c)(y,z)  
  & + \int_\Upsilon k^2 c^{-2}(x) G(x,y) \tildadj(x,y,z)
                 - \nabla_x R(x,y) \cdot \nabla_x \tildadj(x,y,z) \:\dd x \\
  & + \int_{\Gamma_2} \ii k_0 R(x,y) \tildadj(x,y,z)\:\dd \mu(x).
\end{aligned}
\end{equation}
Here, $\tildadj$ has the role of a Lagrange multiplier and 
a specific choice of it, $\adjoint$, will be specified later.
By letting $G=\Gc$ be the solution of the forward problem, 
and hence $R=\Rc$, we can form the Fr\'echet derivative. 
For the sake of brevity, we use the following notation 
\begin{equation}
  G' = \dGc[c]\dc,
\end{equation}
and we omit the variables $x,y,z$ (keeping in mind that 
$\Gobs$ depends on $(x,z)$, $\Gc$ and $G'$ on 
$(x,y)$, and $\tildadj$ on $(x,y,z)$).

\begin{equation} \label{equation:deriv_lagrangian}
\begin{aligned}
D_c \misfitI[c] \dc 
   & = \Ree \big( D_c \lagrangian(c,G,\tildadj) \dc \big)\big|_{G=\Gc}\\
   & = \Ree \bigg( 2\Subar
             \int_{\Sigma} \big(G' \partial_{\nu_x} \Gobs - 
                           \Gobs \partial_{\nu_x}G' \big) \:\dd\mu(x) \\
   & \hspace*{1cm} + \int_\Upsilon \Big(  k^2 c^{-2} G' \tildadj
                         -2k^2 c^{-3} \Gc \tildadj \dc 
                         -\nabla_x G' \cdot \nabla_x \tildadj \Big) \:\dd x \\
   & \hspace*{1cm} + \int_{\Gamma_2} \ii k_0 G' \tildadj \:\dd \mu(x) \bigg).
\end{aligned}
\end{equation}
Grouping together all the terms containing $G'$ and 
replacing it by an arbitrary test function 
$v \in H^1_{\Gamma_1}(\Upsilon)$, the adjoint state
$\adjoint$ is chosen as the solution to
\begin{equation} \label{equation:adjoint_state1}
\begin{aligned}
    2\Subar \int_{\Sigma} \big(v \partial_{\nu_x} \Gobs - 
                          \Gobs \partial_{\nu_x}v \big) \:\dd \mu(x) 
          & + \int_\Upsilon \Big(  k^2 c^{-2} \vbar  \adjoint
                          -\nabla_x \vbar \cdot \nabla_x \adjoint \Big)  \:\dd x \\
          & + \int_{\Gamma_2} \ii k_0 \vbar \adjoint \:\dd \mu(x) = 0.
\end{aligned}
\end{equation}
Note that the first term,
\begin{equation}
v \rightarrow 2\Subar \int_{\Sigma} \big(v \partial_{\nu_x} \Gobs - 
                \Gobs \partial_{\nu_x}v \big) \:\dd \mu(x),
\end{equation}
is a bounded linear functional of $H^1_{\Gamma_1}(\Upsilon)$.
Hence, by the arguments already mentioned in 
\cite[Proposition 3.1]{AdHGS}, there exists a unique 
solution $\adjoint\in H^1_{\Gamma_1}(\Upsilon)$ to 
problem~\eqref{equation:adjoint_state1}. \\

With this choice of adjoint state, \eqref{equation:deriv_lagrangian} 
reduces to
\begin{equation}
 D_c \misfitI[c] \dc = 
                       \Ree\bigg(\int_\Upsilon -2k^2 c^{-3} \Gc \adjoint \dc \:\dd x
                            \bigg).
\end{equation}
Reassembling with $(y,z) \in K_1 \times K_1$, 
we get
\begin{equation}
 \nabla_c \misfit(x)= - \Ree \bigg( \int_{K_1 \times K_1} 
                       2k^2 c^{-3}(x)\Gc(x,y) 
                                   \adjoint(x,y,z) \:\dd y\:\dd z \bigg).
\end{equation}
Note that $\Gc$ is independent of $z$, hence, 
posing
\begin{equation}
  \widehat{\adjoint}(x,y) = \int_{K_1} \adjoint(x,y,z)\:\dd z,
\end{equation}
we have that $\widehat{\adjoint}$ verifies, 
for every $v \in H^1_{\Gamma_1}(\Upsilon)$, 
\begin{equation} \label{equation:adjoint_state_full}
\begin{aligned}
    2 \int_{K_1} \Subar & \int_{\Sigma} \big(v \partial_{\nu_x} \Gobs - 
                        \Gobs \partial_{\nu_x}v \big) \:\dd \mu(x) \:\dd z \\
           + & \int_\Upsilon \Big(  k^2 c^{-2} \vbar  \widehat{\adjoint}
                          -\nabla_x \vbar \cdot \nabla_x \widehat{\adjoint} \Big)  \:\dd x
           + \int_{\Gamma_2} \ii k_0 \vbar \widehat{\adjoint} \:\dd \mu(x) = 0.
\end{aligned}
\end{equation}


\section{Computational experiments}
\label{section:numerical}
\begin{figure}[h!]
\centering
\includegraphics[scale=1]{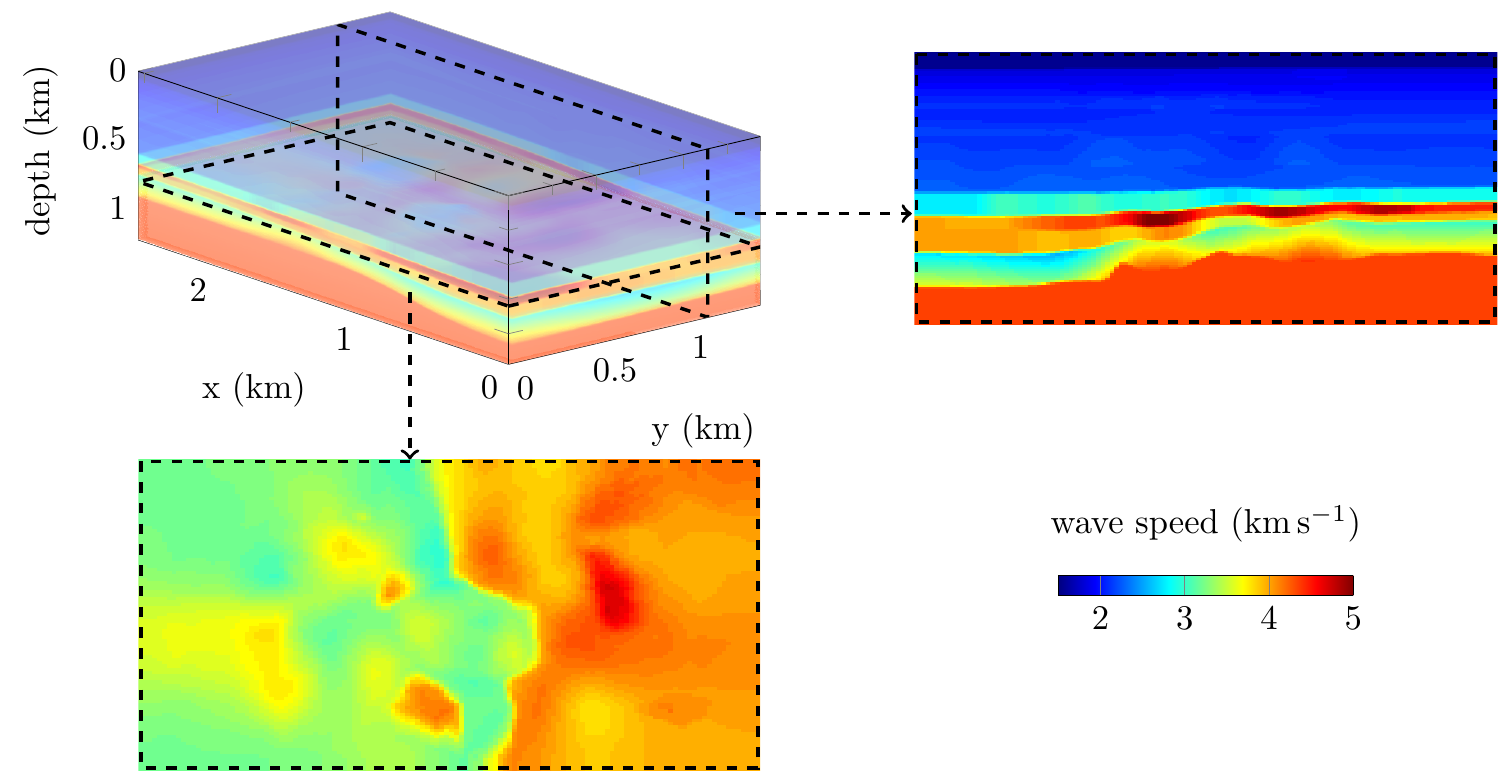}
\caption{Three-dimensional representation, horizontal section at
  $800$\si{\meter} depth and vertical section at $y=1.125$\si{\km} of
  the reference wave speed. It is represented with \num{1527168}
  nodal values.}
\label{fig:statoil}
\end{figure}

For the computational experiments, 
the space coordinates will be denoted by $(x,y,z)$ 
instead of $(x_1,x_2,x_3)$. Let us emphasize that
the $z$ coordinate is conventionally (in geophysical setup)
oriented downwards and can be seen as the depth
of the medium.
We first consider a three-dimensional model
(courtesy Statoil), which is illustrated in
Figure~\ref{fig:statoil}. To have a clear visualization of the
wave speed structures, we show horizontal and vertical sections at
$z=800$\si{\meter} and $y=1.125$\si{\km} respectively. The model is of
size $2.55 \times 1.45 \times 1.22$\si{\km} with variations of wave 
speed from $1500$ to $5200$\si{\meter\per\second}.  We assume that
the density is constant with $\rho=1000$\si{\kilo\gram\per\meter\cubed}.

\begin{figure}[h!] \centering
\subfigure[Three-dimensional computational domain.]
          {\includegraphics[scale=.95]
          {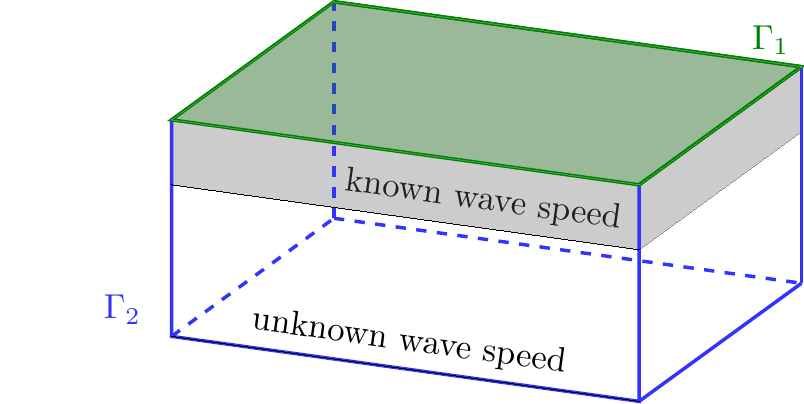}} \hfill
\subfigure[Acquisition with sources contained in a three-dimensional region.]
          {\includegraphics[scale=.95]
          {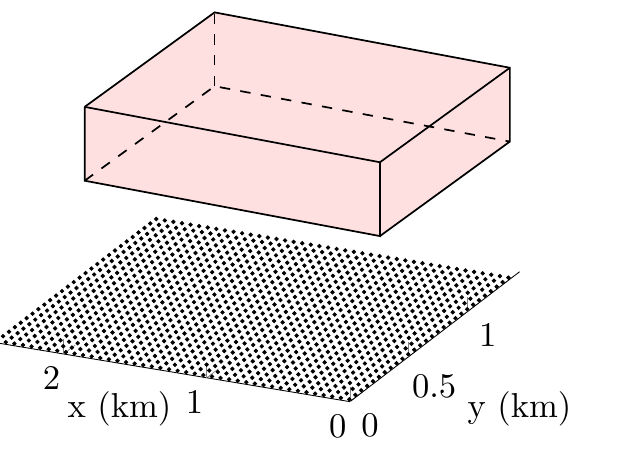}
           \label{fig:setup_sketch_b}} \\[0mm]
\subfigure[Acquisition with sources contained on a plane.]
          {\includegraphics[scale=.95]
          {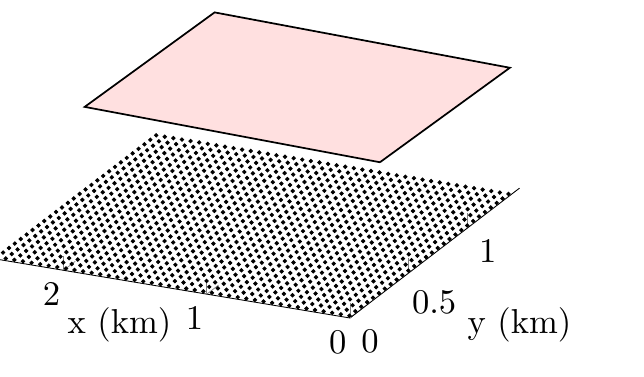}
           \label{fig:setup_sketch_c}}
\caption{Illustration of the configuration. 
         (a) We apply a Dirichlet
         boundary on the upper surface $\Gamma_1$ 
         (in green) 
         and absorbing boundary conditions on 
         $\Gamma_2 = \partial\Upsilon \backslash \overline\Gamma_1$
         (indicated in blue). 
         (b)--(c) The sources that probe the Cauchy data lie 
         in between the receivers and the free surface. 
         Sources are positioned in a three or two-dimensional 
         region (in red). The position of the receivers 
         (black dots) remains fixed. Both receivers and sources lie 
         in the area of known wave speed.}
\label{fig:setup_sketch}
\end{figure}

The seismic acquisition consists of \si{\num{160}} sources and
\si{\num{1376}} receivers, that is, dual sensors recording the Cauchy
data probed by the sources. The receivers are positioned on a regular
$43$ (along the $x$-axis) by $32$ (along the $y$-axis) grid at a fixed
depth below the sources lattice. The configuration is illustrated in
Figure~\ref{fig:setup_sketch}. 
We consider two situations for the discretized sources map:
the sources are first contained in a solid region, in accordance 
with the above analysis (see Figure~\ref{fig:setup_sketch_b}). 
Then, they are restricted on a two-dimensional lattice
(see Figure~\ref{fig:setup_sketch_c}). The first approach,
less common in seismic applications, is yet possible with recent
acquisition technique described in Footnote~\ref{footnote_topseis},
for which the depth of the sources can vary.
Following the situation 
prescribed in Subsection~\ref{subsection:cauchy_data_misfit},
we assume that the uppermost part of the model (which is water), 
in which the Cauchy data are obtained, is known prior to the reconstruction. 
However, we do not assume the knowledge of the wave speed onto the lateral
and bottom boundaries.

We impose a free surface, Dirichlet boundary condition on the
top part, $\Gamma_1$, of $\partial\Upsilon$ and absorbing boundary
conditions on $\Gamma_2 = \partial\Upsilon \backslash \overline\Gamma_1$,
by taking $k_0 =  k c_0^{-1}$ in the third equation of 
\eqref{Green_third_kind_statement}, following 
Engquist and Majda \cite{Engquist1977}.

\begin{figure}[ht!] \centering
  {\includegraphics[scale=1]{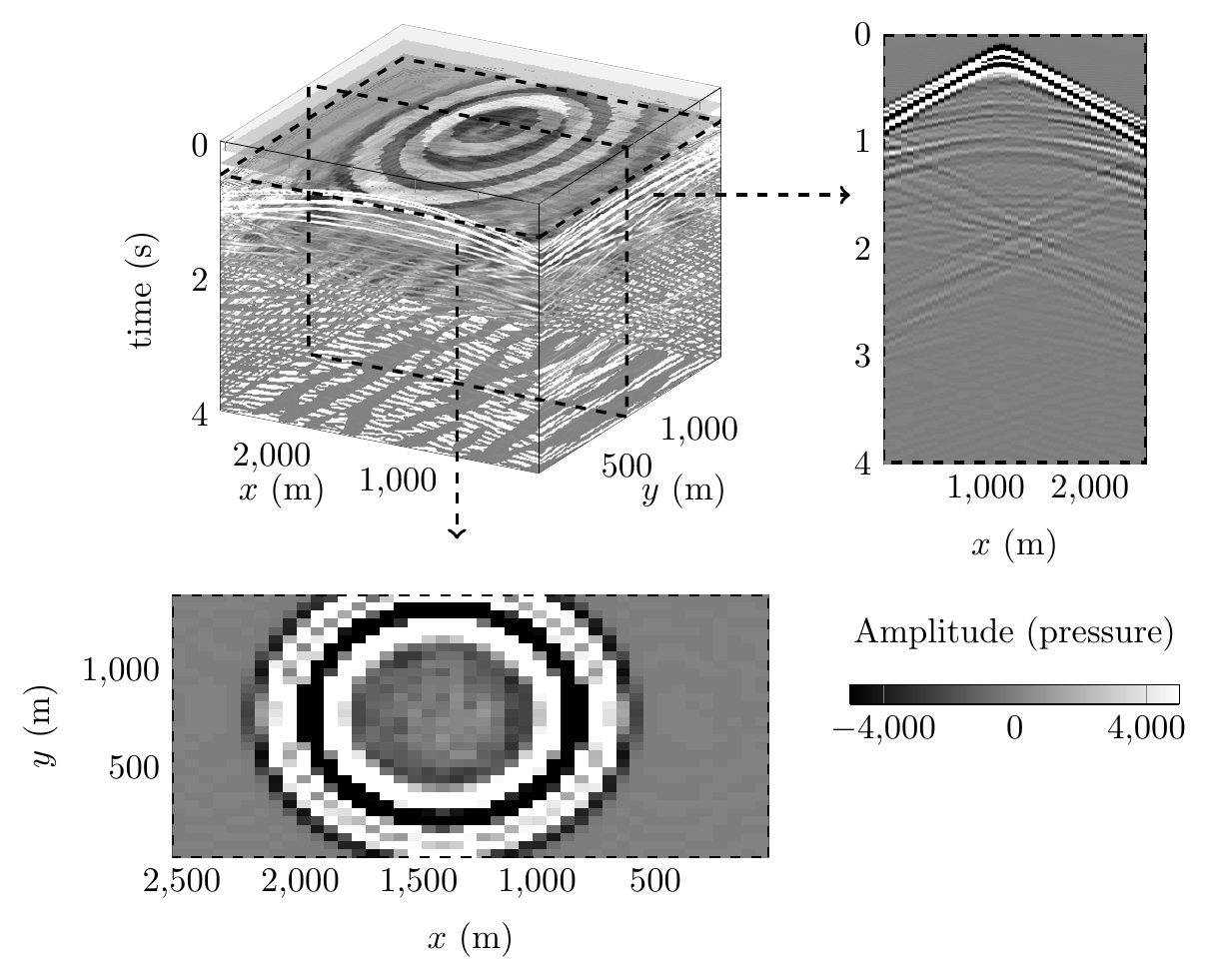}}
\caption{Three-dimensional time-domain pressure trace for a centrally located 
         source. Sections at fixed time $t=0.5$\si{\second} and for 
         a fixed line of receivers positioned in $y=695$\si{\meter}
         are respectively given at the bottom and right of the $3$D 
         visualization. The $x$ and $y$ axis correspond with the 
         receivers map (\textit{i.e.}, every $60$\si{\meter} for $x$ and $45$\si{\meter}
         for $y$).}
\label{fig:td_data1}
\end{figure}

\begin{figure}[ht!] \centering
  {\includegraphics[scale=1]{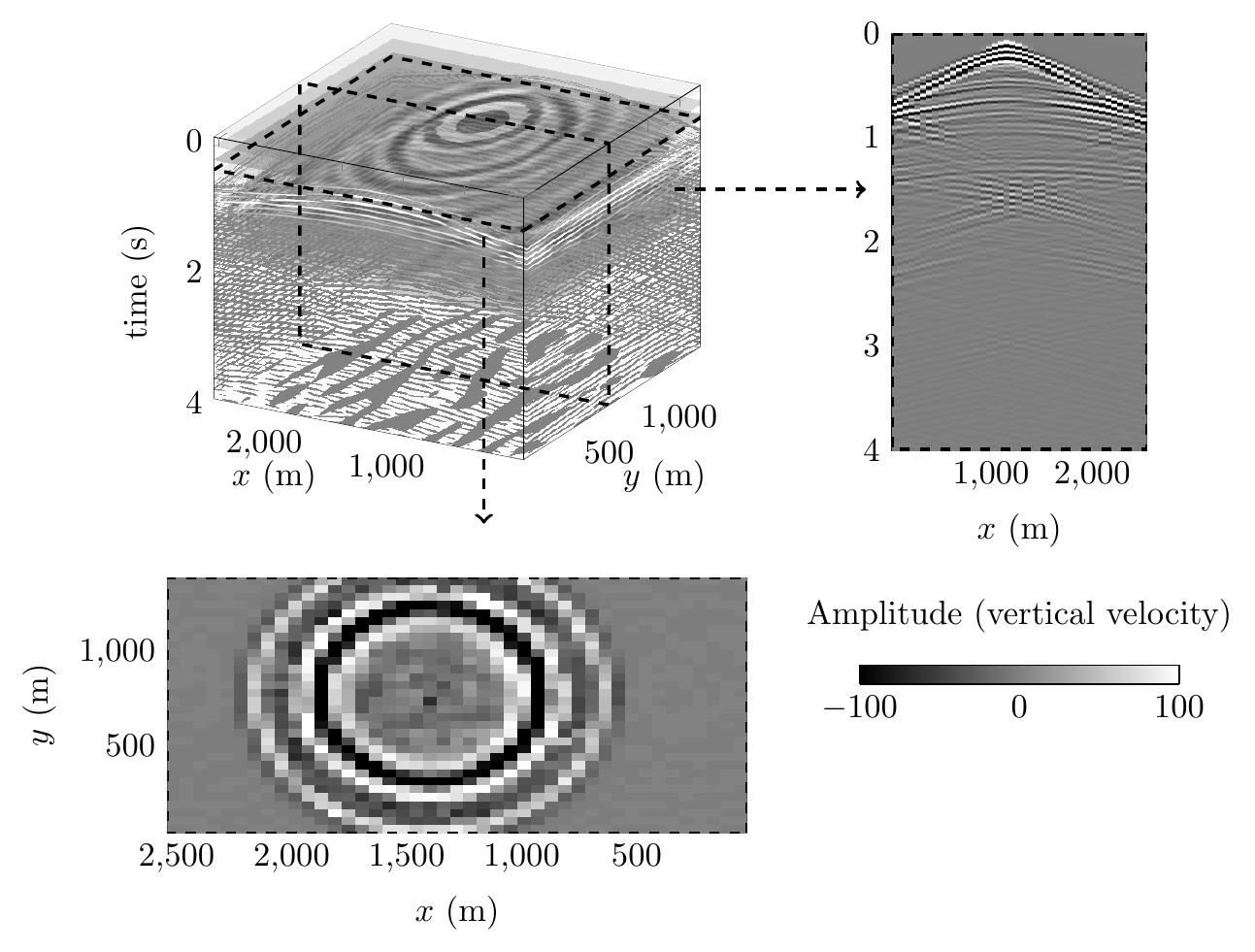}}
\caption{Three-dimensional time-domain vertical velocity trace for a centrally located 
         source. Sections at fixed time $t=0.5$\si{\second} and for 
         a fixed line of receivers positioned in $y=695$\si{\meter}
         are respectively given at the bottom and right of the $3$D 
         visualization. The $x$ and $y$ axis correspond with the 
         receivers map (\textit{i.e.}, every $60$\si{\meter} for $x$ and $45$\si{\meter}
         for $y$).}
\label{fig:td_data2}
\end{figure}

Synthetic dual-sensor data are generated in the time-domain using a
Discontinuous Galerkin (DG) finite element method\footnote{The code that
  was used, here, can be found at
  \url{https://team.inria.fr/magique3d/software/hou10ni/}.}. 
The original data for a single centered source are presented in
Figures~\ref{fig:td_data1} and~\ref{fig:td_data2} for the pressure
and vertical velocity respectively. In the figures, 
we can observe the difference of scale in the amplitudes 
between the pressure and the vertical velocity. 

\setlength{\modelwidth}{8.5cm}
\graphicspath{{figures/statoil/time-domain_data/}}

\setlength{\modelwidth}{8.5cm}
\graphicspath{{figures/statoil/time-domain_data/}}
\begin{figure}[ht!] \centering
 \renewcommand{\modelfile}{trace_scale5e3_statoil-v4000_y16}
 \subfigure[Noiseless synthetic trace section.]
           {\includegraphics[scale=1]{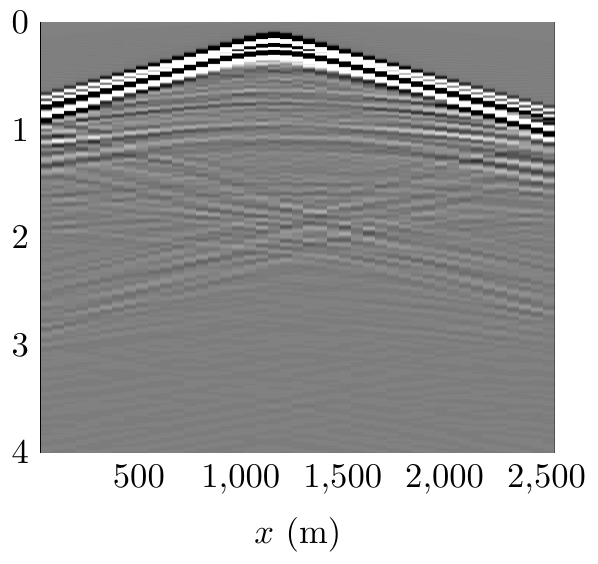}}
 \renewcommand{\modelfile}{trace_scale5e3_statoil-v4000_y16_noise}
 \subfigure[Noisy synthetic trace section using 
            $15$ \si{\decibel} signal to noise ratio.
            It is used to generate the frequency data
            for the reconstruction.]
           {\includegraphics[scale=1]{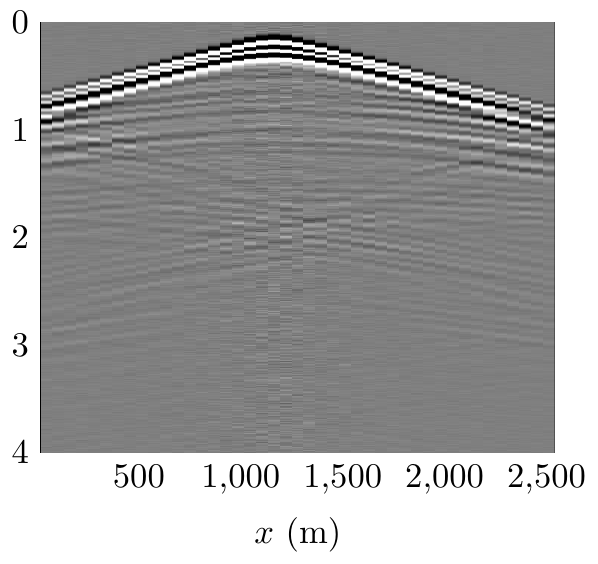}}
\caption{Comparison of noiseless and noisy data for a section 
         of the pressure trace given Figure~\ref{fig:td_data1}.
         The noise is independently generated
         for every receiver and every
         source in the acquisition.}
 \label{fig:td_data1_noise}
\end{figure}

We subject the (time-domain) data to 
Gaussian white noise, using a signal-to-noise ratio 
of $15$\si{\decibel} (this process
is illustrated Figure~\ref{fig:td_data1_noise}).
Note that every receiver for each source
has an independent white noise signal added.
We apply the Fourier transform to these
noisy data and obtain time-harmonic data. 
The effect of noise affects particularly the 
low-frequency regime in seismic, and frequencies
below $3$ \si{\Hz} are usually unusable.

\begin{figure}[ht!] \centering
  \includegraphics[scale=.40]{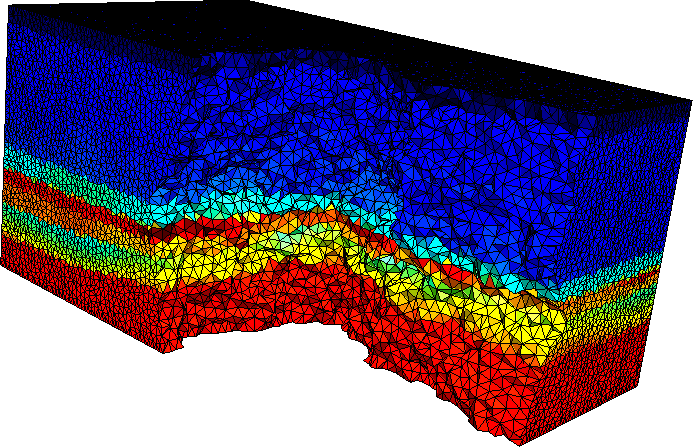}
  \caption{Illustration of the unstructured tetrahedral mesh 
           of the model, using \num{433979} tetrahedra.}
\label{fig:fem_3Dmesh_velocity}
\end{figure}

In this experiment we only select 
$10$\si{\Hz} frequency data
for the reconstruction algorithm and
underlying iterative minimization of the misfit 
functional~\eqref{misfit}.
We simulate time-harmonic data using a Continuous Galerkin 
finite element method (CG). We use an approach and implementation 
similar to the one discussed in Shi \textit{et al.}
\cite{Shi-2017}. The relevant system of equations are solved with the
direct structured solver \textsc{Mumps}, \cite{Amestoy2006}. 
The numerical discretization introduces a tetrahedral 
representation of the model, which we illustrate in 
Figure~\ref{fig:fem_3Dmesh_velocity}.

The choice of CG
(instead of DG) is motivated by 
the memory cost of solving large linear system, 
which is a specificity of the harmonic case
(DG is used for the 
time-domain discretization). In our experiments,
we use order $3$ polynomials to guarantee the 
accuracy (by taking at least four degrees of 
freedom per wavelength, according to the 
common heuristic). Note also that the mesh employed
to generate the synthetic (time-domain) data differs 
from the one used for the inverse (harmonic) problem:
the one to generate the data is refined to 
make sure we consider acutely the variations of the 
reference wave speed model.

For the reconstruction, the wave speed 
uses a piecewise linear representation, 
following~\eqref{potential 2} and \cite{AdHGS}.
For the construction of the hierarchy of stable subspaces, 
the domain partition determining the piecewise linear 
representation of the wave speed is typically significantly coarser
than the tetrahedral mesh. Hence, we define every subdomain, $D_j$, as
the union of mesh elements, $K_i$, according to

\begin{equation}
   D_j = \bigcup_{i=1}^{N_j} K_i,
\end{equation}
where $N_j$ denotes the number of mesh elements in $D_j$.

\begin{figure}[h!]
\centering
\renewcommand{\modelfile}{cp_compression}
{\includegraphics[scale=1]{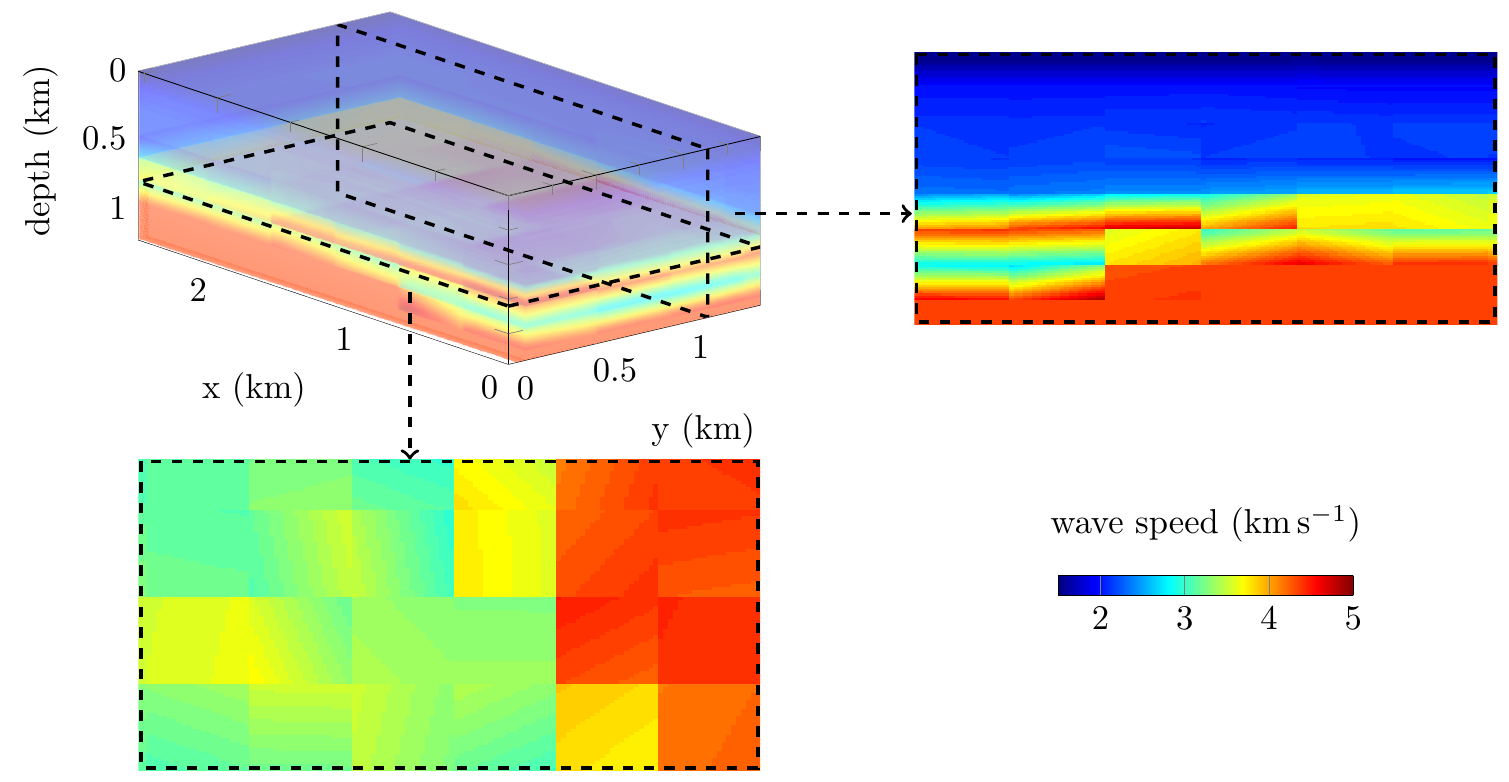}}
\vspace*{0mm}
\caption{Illustration of three-dimensional partitioning of the wave speed
         model of Figure~\ref{fig:statoil}. 
         The size of subdomains is limited to $400$\si{\meter}
         in the $x$ and $y$ directions and $150$\si{\meter}
         in the $z$ direction. This leads to a decomposition of 
         $N=\num{224}$ subdomains where piecewise linear functions
         are used to represent the wave speed.
         }
\label{fig:statoil_compression}
\end{figure}

To achieve the decomposition, we apply a structured 
decision where the maximal size of the subdomains is
chosen in every direction to define the subspace. 
Then, piecewise linear functions are 
employed onto each generated subdomain to represent the model. We 
illustrate the effect of piecewise linear partitioning 
applied on the wave speed model in Figure~\ref{fig:statoil_compression}, 
where the size of the subdomain is at most $400$\si{\meter}
in the $x$ and $y$ directions, and $150$\si{\meter}
in the $z$ direction; this leads to a decomposition 
with $N=\num{224}$ subdomains 
and \num{896} coefficients to represent the model.
Inherent model error is introduced from those two levels of 
representation (the mesh and the partitioning).
Because we do not know the subsurface geometry 
a priori, the piecewise linear partition relies on the gradient 
of the misfit functional instead of the wave speed. Naturally, the 
more subdomains are taken, the more accurate can the representation be. 

We proceed using single-frequency, $10$\si{\Hz} data 
and a fixed domain partition. 
Exploiting the Lipschitz stability result obtained 
in Theorem~\ref{teorema principale} above,
the Landweber iteration \cite{deHoop2012}
provides a convergence analysis. The initial model needs to be within
the radius of convergence. 
Algorithm~\ref{algo:fwi} summarizes the steps of the 
reconstruction procedure. 

\begin{algorithm}[ht!] \setstretch{1.25} 
\SetKwInput{KwData}  {Preliminary material}
\SetKwInput{KwResult}{Results}

\KwData{
\vspace*{-.2\baselineskip}
\begin{itemize} \setlength{\itemsep}{-2pt}
 \item time-domain observation data at the receivers location,
 \item user prescribed frequency $k$ 
       and associated partition $N$,
 \item user prescribed initial model $c_1$ and number of 
       iterations $n_\text{iter}^{(\min)}$, $n_\text{iter}^{(\max)}$,
 \item user prescribed stagnation parameters $\epsilon_\misfit$ 
       and $n_{\epsilon}$.
\end{itemize}
}
Computation of the Fourier transform of the time-domain Cauchy traces at $k$. \\
\BlankLine 
  \textbf{Optimization loop} \For{$j \in \{1, \ldots, n_\text{iter}^{(\max)} \}$}{
  \begin{itemize}[leftmargin=*]\setlength{\itemsep}{-2pt}
    \item solve the Helmholtz equation~\eqref{Green_third_kind_statement} 
          at frequency $k$ with wave speed $c_j$; \\
    \item compute the misfit functional~\eqref{misfit} from 
        the simulation and observation data;      \\
    \item compute the gradient of the misfit functional, 
        $\nabla\misfit(c_j)$ with the adjoint-state method (see 
        Subsection~\ref{subsec:adjoint_state_misfit}); \\
    \item compute the search direction, $s_j$, here
        we use the nonlinear conjugate gradient method with 
        Polak--Ribi\`ere formula (\textit{cf.} \cite[Section~5.2]{Nocedal2006}); \\
    \item apply segmentation onto the search direction 
         (see illustration on Figure~\ref{fig:statoil_compression}); \\ 
    \item compute the step length $\alpha$ with line search method 
          (backtracking, \textit{cf.} \cite[Chapter~3]{Nocedal2006}); \\
    \item model update $c_{j+1} = c_j - \alpha s_j$; \\
  \end{itemize}
  \vspace*{-.5\baselineskip}
  \If{$j \geq n_\text{iter}^{(\min)}$ \textnormal{\textbf{and}} $j > n_{\epsilon}$ }{
    \begin{itemize}[leftmargin=*]\setlength{\itemsep}{-2pt}
      \item compute the stagnation criterion
       $\mathfrak{e} = \dfrac{\misfit(c_{j-n_{\epsilon}})
             -\misfit(c_{j})}{\misfit(c_{j-n_{\epsilon}})};
       $              
      \item \textnormal{\textbf{if}} 
            $(\mathfrak{e} < \epsilon_\misfit)$: \textbf{exit optimization loop}.
    \end{itemize}
   }
  }

 \BlankLine 
 \caption{algorithm for the reconstruction of subsurface parameters using 
          piecewise linear model partition and Cauchy data.}
 \label{algo:fwi}
\end{algorithm}

Two parameters can decide of the termination of the
procedure: if the number of iterations $n_\text{iter}^{(\max)}$ is
reached, or if the cost function stagnates (criteria 
$\epsilon_\misfit$ and $n_{\epsilon}$), 
see Algorithm~\ref{algo:fwi}. In the following
experiments, we impose
\begin{equation}
\begin{aligned}
      n_\text{iter}^{(\min)} = 50, \quad n_\text{iter}^{(\max)} = 250, \\
\quad n_{\epsilon}=10 \quad \text{and} \quad 
\epsilon_\misfit=0.01.
\end{aligned}
\end{equation}
Hence, the number of iterations is kept relatively high
and the stagnation stops the procedure. 
More precisely with the given numbers, the algorithm stops if the 
difference in the misfit functional over the last ten iterations 
is less than $1$\%.

\subsection{Single-frequency data}

\begin{figure}[ht!] \centering
\renewcommand{\modelfile}{vp_start_smooth}
{\includegraphics[scale=1]{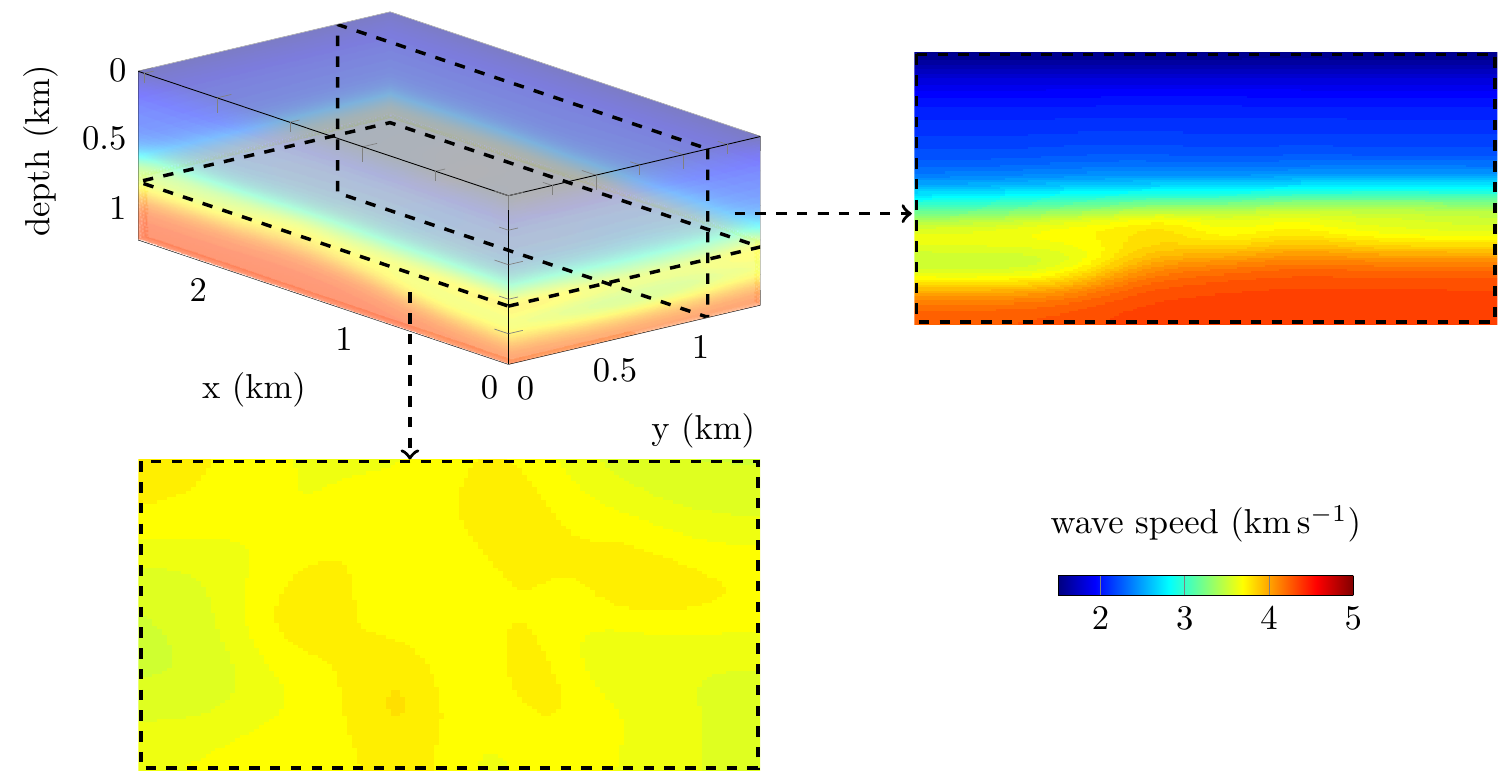}}
\caption{Three-dimensional representation, horizontal section at
  $800$\si{\meter} depth and vertical section at $y=1.125$\si{\km} of
  the initial wave speed model.}
\label{fig:statoil_start_smooth}
\end{figure}

From the Cauchy data at $10$\si{\Hz}, we carry out a reconstruction of
the reference model starting from the smooth model depicted
Figure~\ref{fig:statoil_start_smooth}. We encode the principal
variation and appropriate order of magnitude of the wave speed in the
initial model.

\begin{figure}[ht!]
\centering
\subfigure[Reconstruction where the sources in 
           the acquisition are positioned in a 
           three-dimensional area (see 
           Figure~\ref{fig:setup_sketch_b}); 
           the relative $L^2$ difference 
           with the reference model Figure~\ref{fig:statoil} is 
           $\mathcal{E}_\text{rel}=0.083$.
           ]
          {{\includegraphics[scale=1]{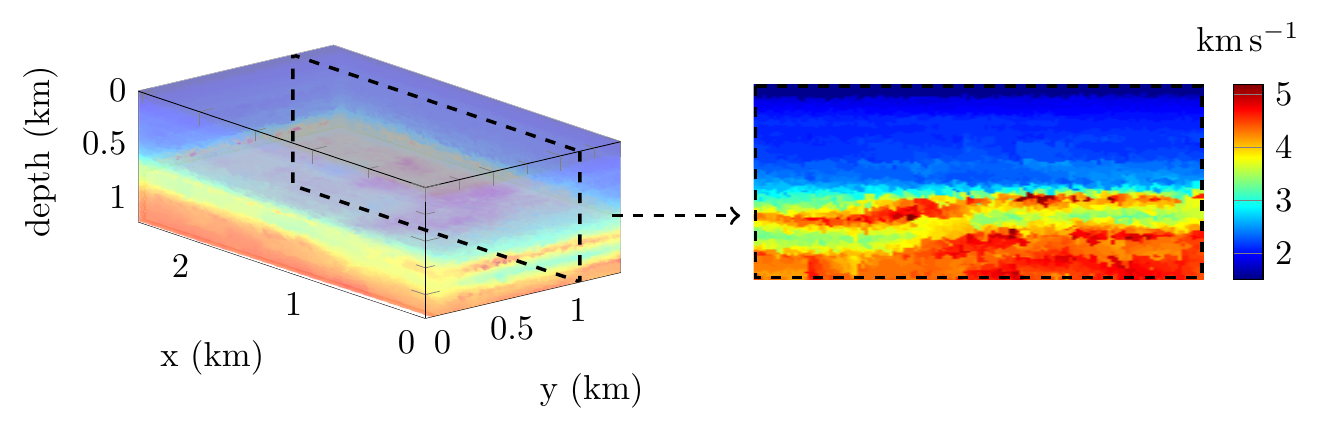}}
           \label{fig:fwi_single-frequency_volume}}\\[-5mm]
\renewcommand{\modelfile}{cp_unifreq_10hz}
\subfigure[Reconstruction where the sources in 
           the acquisition are limited on a 
           plane (see Figure~\ref{fig:setup_sketch_c});
           the relative $L^2$ difference 
           with the reference model Figure~\ref{fig:statoil} is 
           $\mathcal{E}_\text{rel}=0.088$.
           ]
          {{\includegraphics[scale=1]{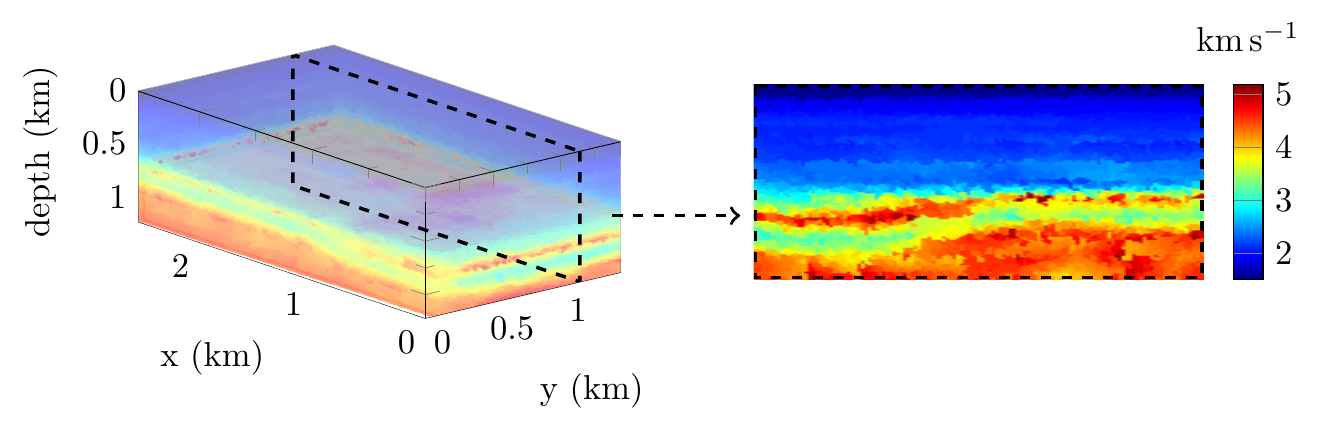}}
           \label{fig:fwi_single-frequency_surface}}
\caption{Three-dimensional representation and
         vertical section at $y=1.125$\si{\km} of
         the reconstruction from $10$\si{\Hz} Cauchy 
         data after $175$ iterations. 
         The partition consists of $N=1089$ subdomains, leading
         to a total number of unknowns of $4\times1089 = 4356$. The initial
         model is shown in Figure~\ref{fig:statoil_start_smooth}.}
\label{fig:fwi_single-frequency}
\end{figure}

The model is partitioned in \num{1089} subdomains where piecewise linear
functions are used to define the wave speed. This leads to a total
number of unknowns of \num{4}$\times$\num{1089} = \num{4356} (while the reference medium
has \num{1527168} nodal values, see Figure~\ref{fig:statoil}). The
key, here, is the low-dimensional subspace used for
regularization. The partition is adapted to the gradient computation
via segmentation. We carry out \num{175} iterations for the two situations.
In Figure~\ref{fig:fwi_single-frequency_volume},
we show the reconstruction when the sources are positioned 
in a volume above the receivers (see Figure~\ref{fig:setup_sketch_b});
in Figure~\ref{fig:fwi_single-frequency_surface}, we show the reconstruction
when the sources are restricted on a two-dimensional plane
(see Figure~\ref{fig:setup_sketch_c}).
To compare the accuracy of the reconstructions, we 
use the relative $L^2$ norm of the difference between the 
reference model and the reconstruction:
\begin{equation}
  \text{relative $L^2$ error} = \mathcal{E}_\text{rel} =
          \dfrac{\Vert c_\dagger - c_\text{r} \Vert}{\Vert c_\dagger \Vert},
\end{equation}
where $c_\dagger$ is the reference model of 
Figure~\ref{fig:statoil} and $c_\text{r}$ the
final reconstruction.

The two acquisitions provide accurate recoveries of the subsurface,
with some improvement when the sources are positioned in a volume:
the width of the increased wave speed layer and the deepest 
values of the wave speed are better retrieved with this type of 
acquisition. However, the reconstruction from sources limited on a
plane is very close.
In the subspace, we have drastically reduced the number of unknowns in
the representation as compared with the original representation,
namely to $0.3\%$. Nonetheless, the reconstruction captures the main
features of the model, including the alternation of high and low
values in the vertical direction and the resolution remains reasonable.

\begin{remark}[Improved visualization with Gaussian filtering]
  \label{rk:gaussian_image}
  The visualization of the reconstruction may suffer 
  from the tetrahedral mesh employed for the numerical
  discretization. It is simple to improve the visualization
  by applying a smoothing filter onto the image. 
  This can be done, for example, with the \texttt{imgaussfilt} 
  function of MATLAB, which applies a Gaussian smoothing filter. 
  In Figure~\ref{fig:fwi_single-frequency_smooth}, 
  we show the resulting image when applied onto the reconstruction 
  of Figure~\ref{fig:fwi_single-frequency_surface}. It allows a 
  better identification of the recovered structures.
  Note that this procedure is done a-posteriori, 
  independently of the reconstruction algorithm, and is effortless.
  \bigskip

  \begin{figure}[h!] \centering
  \renewcommand{\modelfile}{cp_unifreq_10hz_smooth}
  {\includegraphics[scale=1]{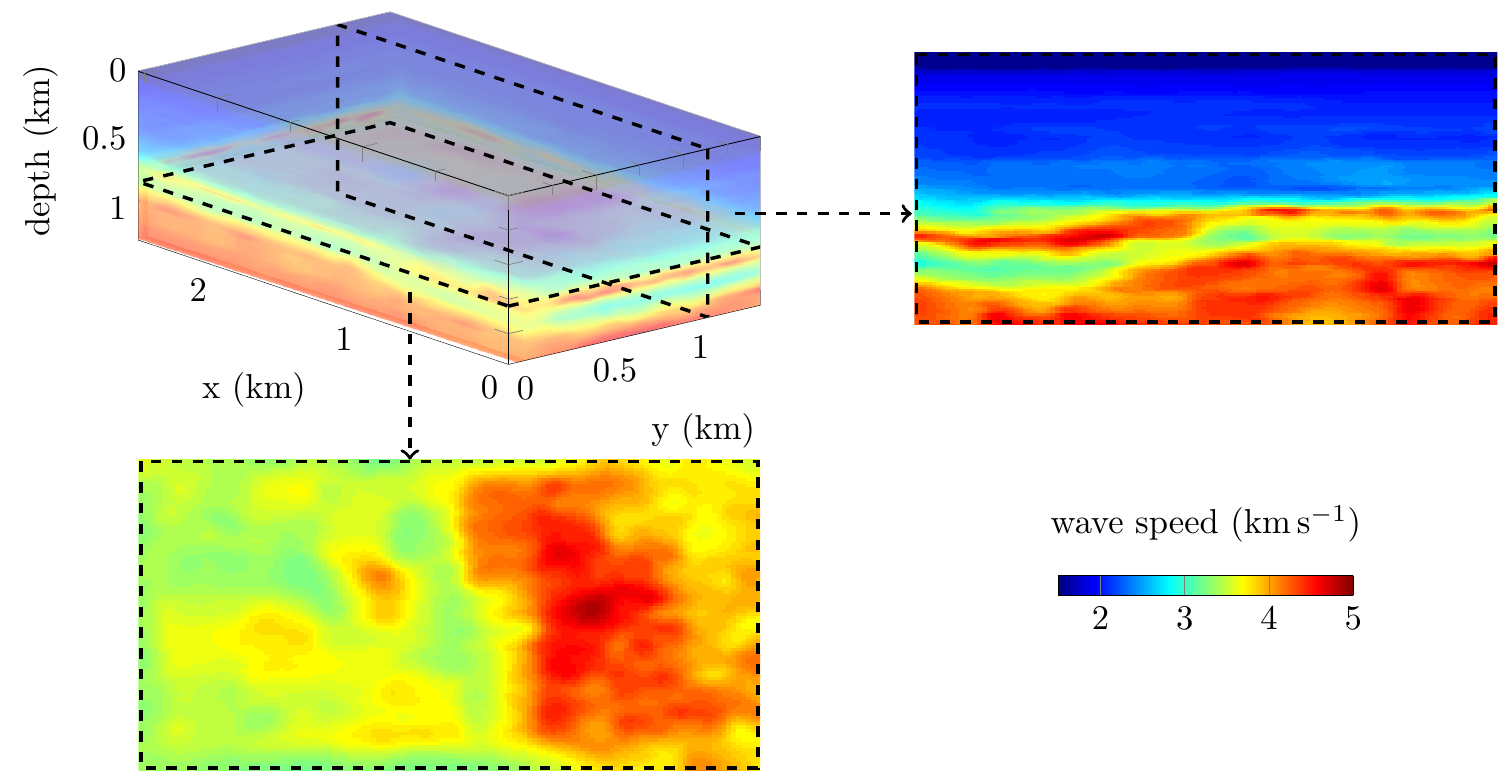}}
  \vspace*{0mm}
  \caption{Gaussian filtering of the reconstruction obtained 
           from $10$\si{\Hz} Cauchy data (Figure~\ref{fig:fwi_single-frequency_surface}),
           three-dimensional representation, horizontal section at
           $800$\si{\meter} depth and vertical section at $y=1.125$\si{\km}.
           The relative $L^2$ difference 
           with the reference model Figure~\ref{fig:statoil} 
           is $\mathcal{E}_\text{rel}=0.079$.
           }
  \label{fig:fwi_single-frequency_smooth}
  \end{figure}
  
\end{remark}

In the following experiments, we only consider the case where
sources are restricted on a plane, following the acquisition
illustrated Figure~\ref{fig:setup_sketch_c}, for simplicity.

\subsection{Single-frequency data, depth varying initial model}

\begin{figure}[ht!]
\centering
\renewcommand{\modelfile}{vp_start_smooth1D}
{\includegraphics[scale=1]{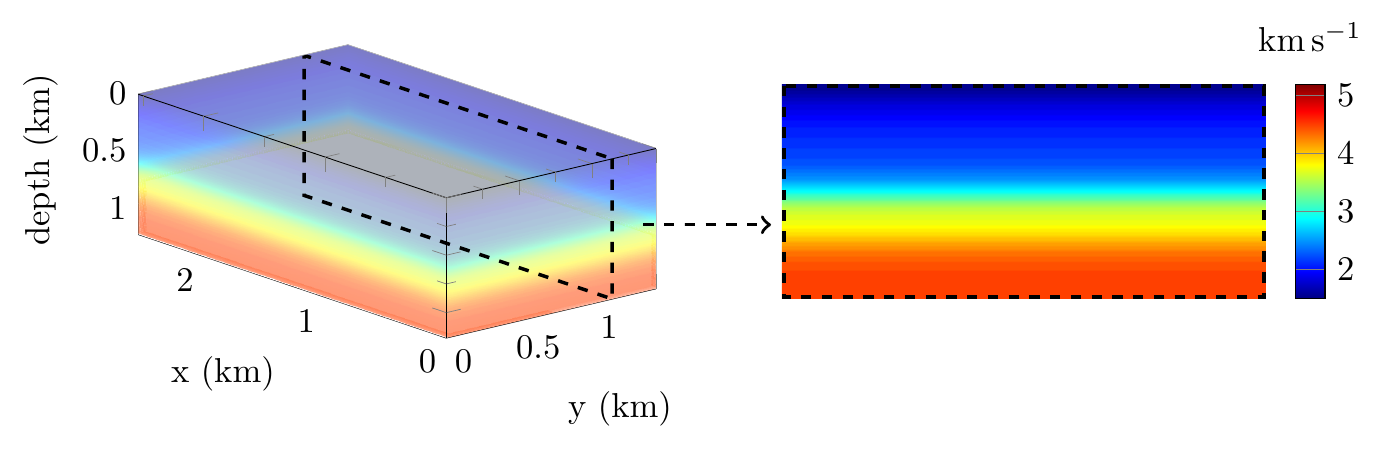}}
\caption{Three-dimensional representation and 
         vertical section at $y=1.125$\si{\km} of
         the initial model, which varies in depth only.}
\label{fig:statoil_start_smooth_1d}
\end{figure}

We repeat the experiment carried out in the previous subsection
(with the two-dimensional sources lattice),
but with a simplified initial model, see 
Figure~\ref{fig:statoil_start_smooth_1d}. That is, 
here, the initial model only contains an indication of the average 
variation of wave speed in depth. The idea behind this experiment 
is to test the radius of convergence on the one hand, and the 
closeness of the true model and the best projection onto a 
low-dimensional stable subspace on the other hand.

\begin{figure}[h!]
\centering
{\includegraphics[scale=1]{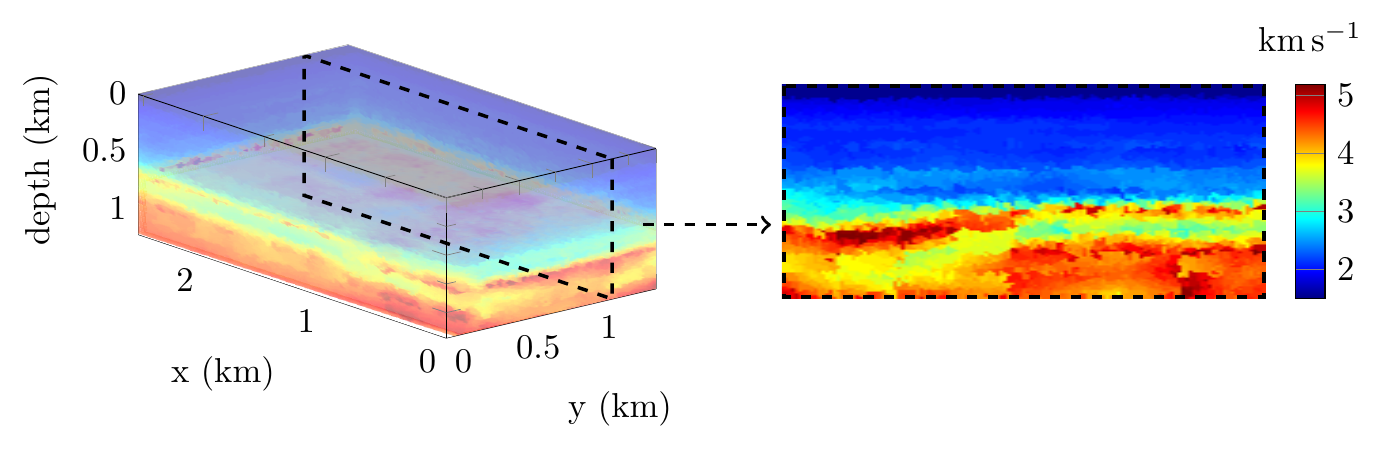}}
\caption{Three-dimensional representation and 
         vertical section at $y=1.125$\si{\km} of
         the reconstruction from $10$\si{\Hz} Cauchy data after $175$
         iterations. The partition consists of $N=1089$ subdomains, leading
         to a total number of unknowns of $4\times1089 = 4356$. The initial
         model varies in depth only and is shown in
         Figure~\ref{fig:statoil_start_smooth_1d}.
         The relative $L^2$ difference 
         with the reference model Figure~\ref{fig:statoil} 
         is $\mathcal{E}_\text{rel}=0.119$.
         }
\label{fig:fwi_single-frequency_1d}
\end{figure}

\begin{figure}[h!]
\centering
{\includegraphics[scale=1]{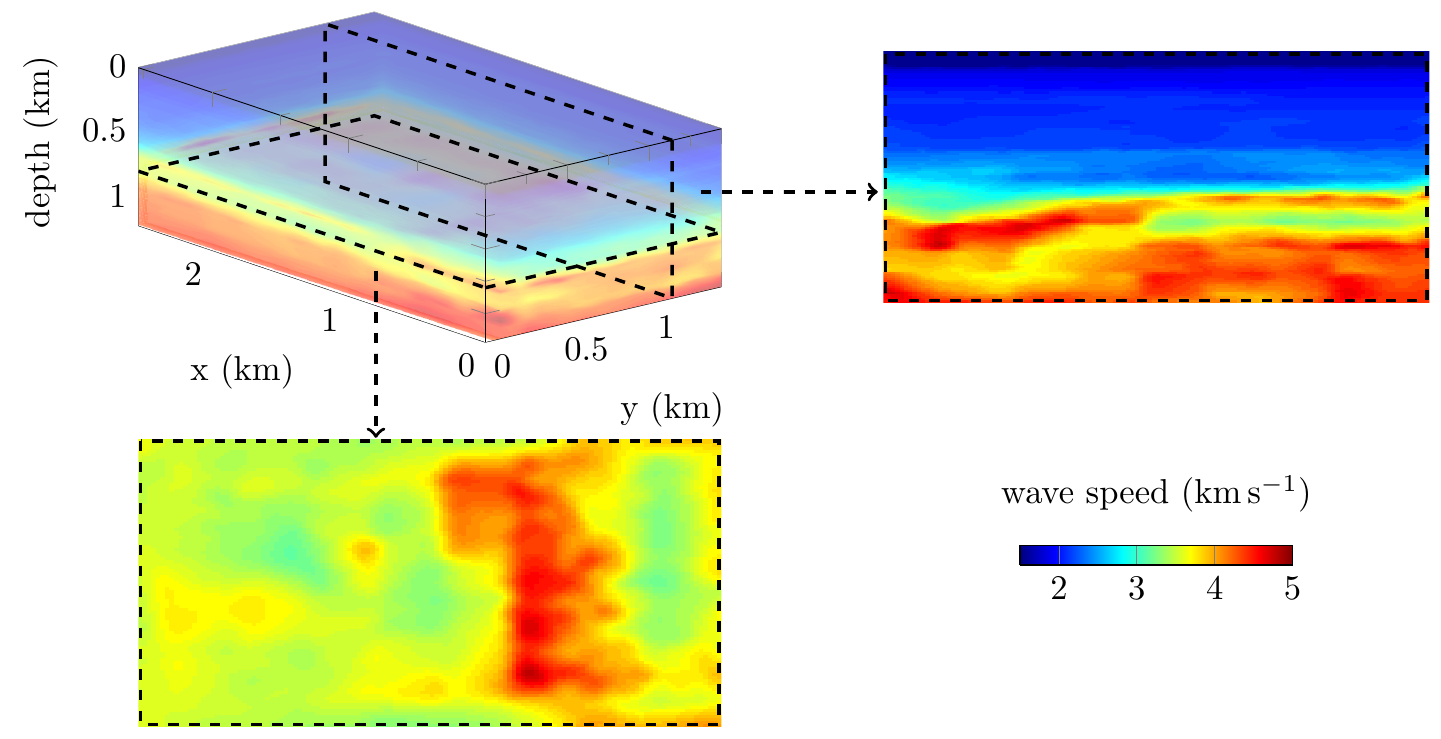}}
\caption{Gaussian filtering (see Remark~\ref{rk:gaussian_image})
         of the reconstruction (Figure~\ref{fig:fwi_single-frequency_1d})
         obtained from $10$\si{\Hz} Cauchy data starting with 
         the initial model that varies in depth only
         (Figure~\ref{fig:statoil_start_smooth_1d}).
         Three-dimensional representation, horizontal 
         section at $800$\si{\meter} depth and vertical 
         section at $y=1.125$\si{\km}.
         The relative $L^2$ difference 
         with the reference model Figure~\ref{fig:statoil} 
         is $\mathcal{E}_\text{rel}=0.108$.
         }
\label{fig:fwi_single-frequency_1d_smooth}
\end{figure}

In Figures~\ref{fig:fwi_single-frequency_1d} 
and~\ref{fig:fwi_single-frequency_1d_smooth}, 
we present the result after $175$ iterations using 
$10$\si{\Hz} Cauchy data. As in the first experiment,
we have \num{1089} subdomains in the partition and
piecewise linear representations. Despite the lack 
of initial information we still retrieve the main 
features and appropriate contrasts in the wave speed. 
However, we lost accuracy as compared with the previous
example (especially on the side), but the deep layer of 
low wave speed is well identified nonetheless. 

\begin{figure}[ht!] \centering
  {\includegraphics[scale=1]{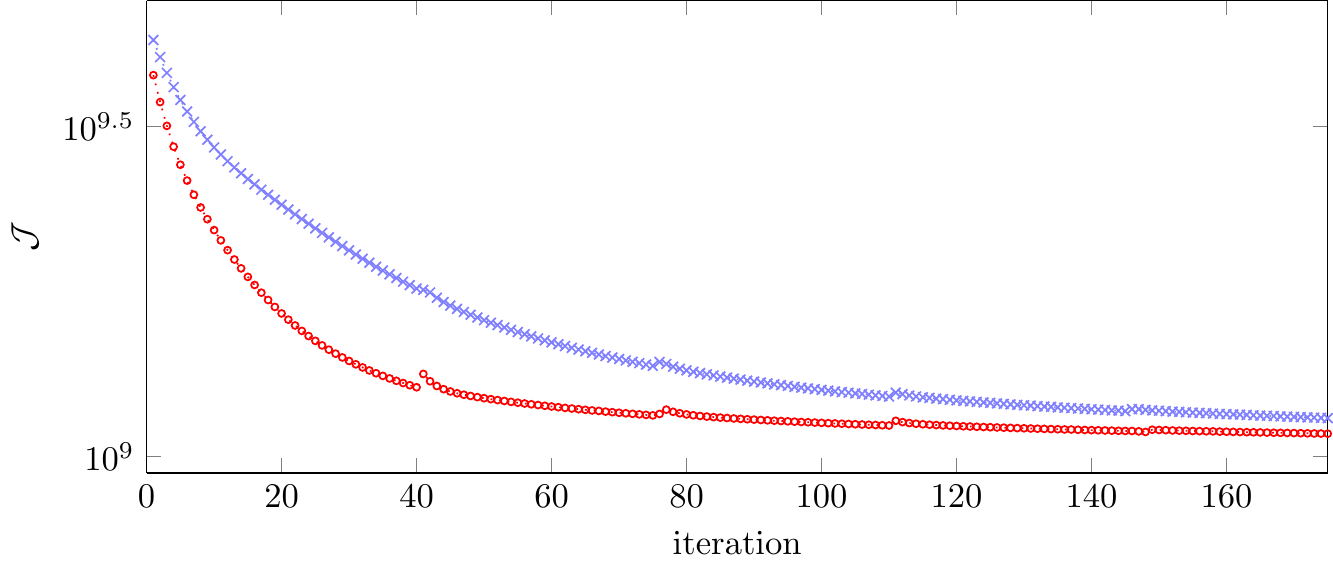}}
  \caption{Evolution of the misfit functional
           with iterations, the blue crosses
           correspond
           with the case where the initial model is
           one-dimensional (see
           Figure~\ref{fig:statoil_start_smooth_1d} 
           for the initial model and 
           Figure~\ref{fig:fwi_single-frequency_1d}
           for the associated reconstruction);
           the red circles 
           correspond with the case where the initial 
           model is smooth (see Figure~\ref{fig:statoil_start_smooth} 
           for the initial model and 
           Figure~\ref{fig:fwi_single-frequency_surface}
           for the associated reconstruction).
           }
  \label{fig:misft_residuals}
\end{figure}

\bigskip
\bigskip
In Figure~\ref{fig:misft_residuals}, we show the 
evolution of the misfit functional with iterations
on a logarithmic scale,
where we compare with the previous experiment
that used a smooth initial model 
(see Figure~\ref{fig:statoil_start_smooth}
and Figure~\ref{fig:fwi_single-frequency_surface} 
for the reconstruction). As expected, the first iteration
using the smooth initial model provides a reduction
in the misfit functional compared to the one-dimensional
starting model. 
The decrease of the misfit functional 
is relatively fast for the initial iterations, especially
when starting with the smooth model, and we observe 
a slow evolution after about $100$ iterations in both 
configurations. Eventually, 
we observe the stagnation which stops the procedure. 
As indicated with the $L^2$ norm of the difference between
the reference model 
(see Figures~\ref{fig:fwi_single-frequency_surface} 
and~\ref{fig:fwi_single-frequency_1d}), starting with 
the smooth model provides a better approximation.

\setlength{\modelwidth}{12cm}

\begin{figure}[ht!] \centering
  \graphicspath{{figures/statoil/data/2D/dirichlet_1e3/}}
  \pgfmathsetmacro{\cmin} {-1000} \pgfmathsetmacro{\cmax} {1000}
  \renewcommand{\modelfile}{data-shot76-10hz-hou10ni} 
  \subfigure[Fourier transform of the time-domain observation.]
  {\scalebox{1}{{\includegraphics[scale=1]{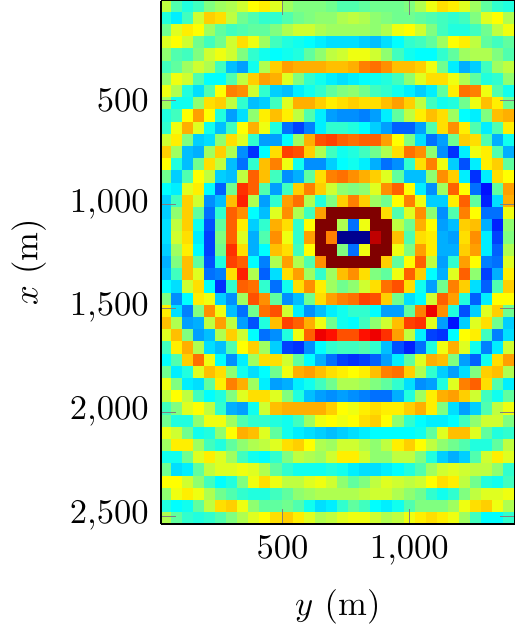}}}} \hfill
  \renewcommand{\modelfile}{data-shot76-10hz-startsmoothWB} 
  \subfigure[Simulated data using the initial 
             model of Figure~\ref{fig:statoil_start_smooth_1d}.]
  {\scalebox{1}{{\includegraphics[scale=1]{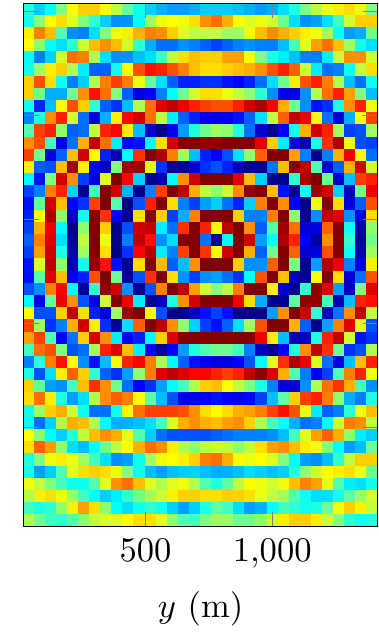}}}}  \hfill
  \renewcommand{\modelfile}{data-shot76-10hz-reconstruction_singlefreq}
  \subfigure[Simulated data from the reconstructed model
             Figure~\ref{fig:fwi_single-frequency_1d}.]
  {\scalebox{1}{{\includegraphics[scale=1]{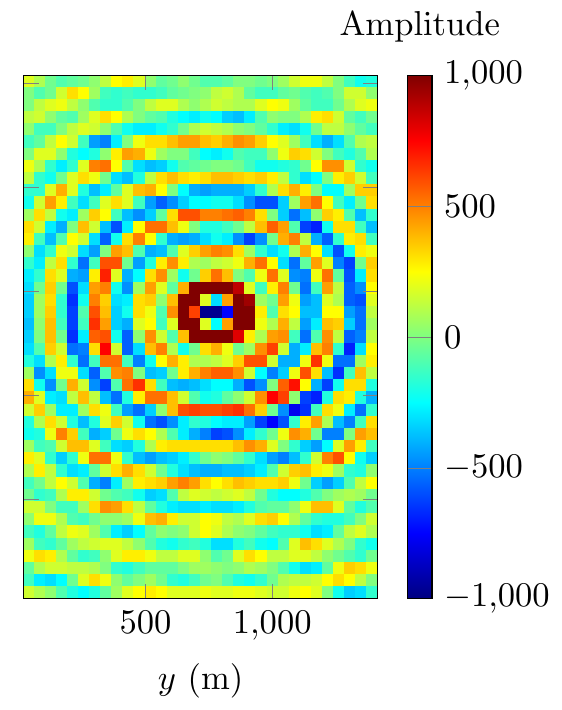}}}}
  \caption{Comparison of the $10$\si{\Hz} frequency pressure 
           data captured at the receivers location for 
           a centrally located source.}
\label{fig:fwi_single-frequency_data_2D_dirichlet}
\end{figure}

\begin{figure}[ht!] \centering
  \graphicspath{{figures/statoil/data/2D/neumann_2e1/}}
  \pgfmathsetmacro{\cmin} {-20} \pgfmathsetmacro{\cmax} {20}
  \renewcommand{\modelfile}{data-shot76-10hz-hou10ni}
  \subfigure[Fourier transform of the time-domain observation.]
  {\scalebox{1}{{\includegraphics[scale=1]{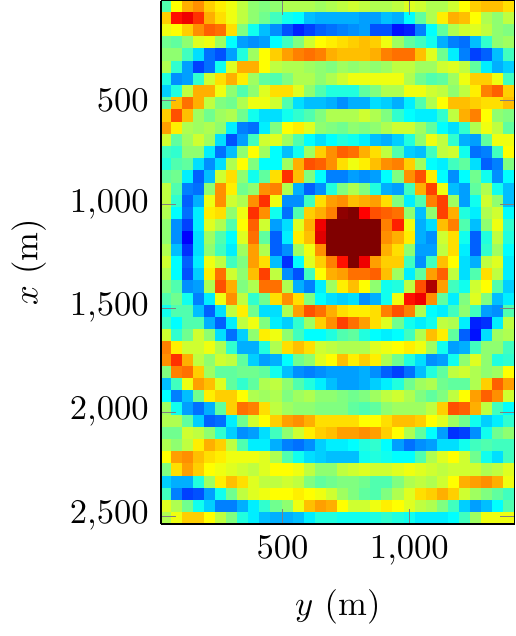}}}} \hfill
  \renewcommand{\modelfile}{data-shot76-10hz-startsmoothWB} 
  \subfigure[Simulated data using the initial 
             model of Figure~\ref{fig:statoil_start_smooth_1d}.]
  {\scalebox{1}{{\includegraphics[scale=1]{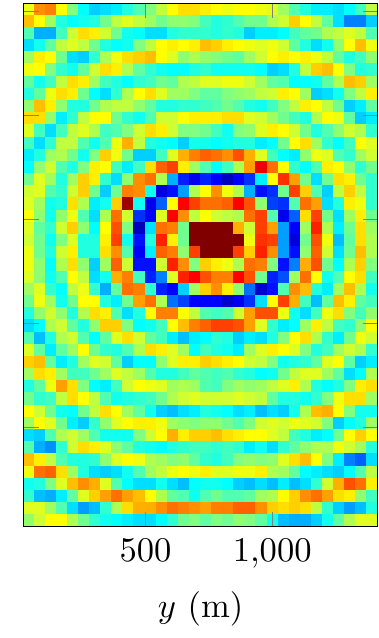}}}}  \hfill
  \renewcommand{\modelfile}{data-shot76-10hz-reconstruction_singlefreq}
  \subfigure[Simulated data from the reconstructed model
             Figure~\ref{fig:fwi_single-frequency_1d}.]
  {\scalebox{1}{{\includegraphics[scale=1]{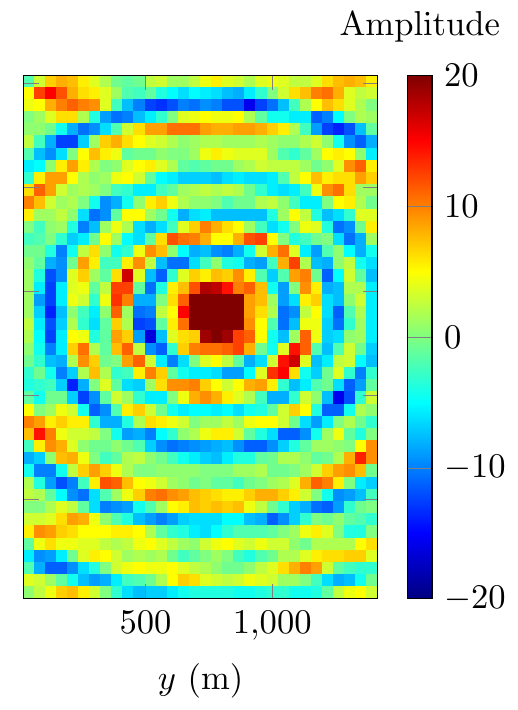}}}}
  \caption{Comparison of the $10$\si{\Hz} frequency vertical 
           velocity data captured at the receivers location 
           for a centrally located source.}
\label{fig:fwi_single-frequency_data_2D_neumann}
\end{figure}

Figures~\ref{fig:fwi_single-frequency_data_2D_dirichlet}
and~\ref{fig:fwi_single-frequency_data_2D_neumann} compare 
the observed, initial and reconstruction data, using the 
full receivers map associated with a centrally located source. 
We see that the data from the recovered wave speed provide 
a pattern that is similar to the Fourier transform of the 
time-domain observations.

\setlength{\modelwidth}{6.90cm}

\section{Perspectives on independent locations of sources in the discretized settings}
\label{section:numerical_independent}
The misfit functional~\eqref{misfit} defined for the Cauchy data 
has an interesting feature, because of the intuitive 
differentiation between acquisition sets for the 
observations and simulations. It is materialized by 
the double integral over $K_1 \times K_1$. 
The perspective is here to separate in, say, $K_2 \times K_1$.

In the usual context of minimization involving the 
direct difference between observations and simulations, 
such as the standard least squares, the setup for simulation is
imposed by the field acquisition (source position and 
wavelet). Consequently, absence of knowledge leads to 
the failure of the algorithm. Here, providing this new
misfit functional, we expect our iterative minimization 
algorithm to be free of those considerations, 
introducing extreme flexibility for the setup, where 
only the position of the receivers is required.

\begin{figure}[h!] 
\centering
\renewcommand{\modelfile}{cp_true}
{\includegraphics[scale=1]{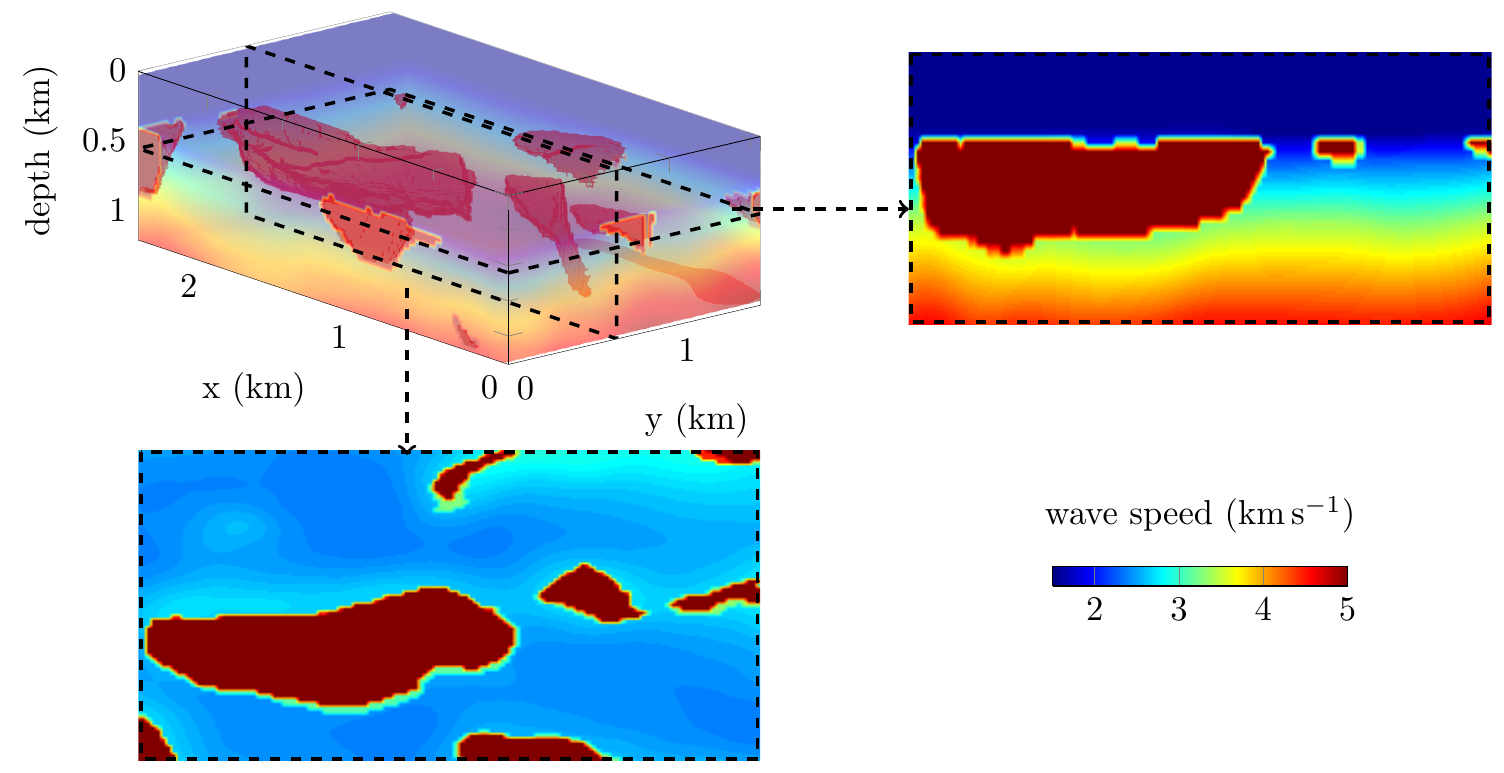}}
\caption{Three-dimensional representation, 
         horizontal section at $550$\si{\meter} 
         depth and vertical section at $y=670$\si{\meter} 
         of the reference wave speed encompassing
         salt domes. The wave speed in the domes is 
         a constant of $4500$\si{\meter\per\second}.}
\label{fig:salt}
\end{figure}

To illustrate the potential of the method,
we design an experiment where the simulation sources
differ from the observational ones. 
We consider a subsurface wave speed model where salt domes
(objects having a large speed contrast) are present. The 
model is illustrated in Figure~\ref{fig:salt} and is of size 
$2.46 \times 1.56 \times 1.2$\si{\km}. It consists
in a smooth background with contrasting objects having a 
wave speed of $4500$\si{\meter\per\second}. We assume that
the density remains constant with 
$\rho=1000$\si{\kilo\gram\per\meter\cubed}. This medium 
is very different in nature from the previous one, and salt 
domes are traditionally challenging in seismic exploration
(especially where their presence is initially unknown).

\begin{figure}[ht!]
\centering
{\includegraphics[scale=1]{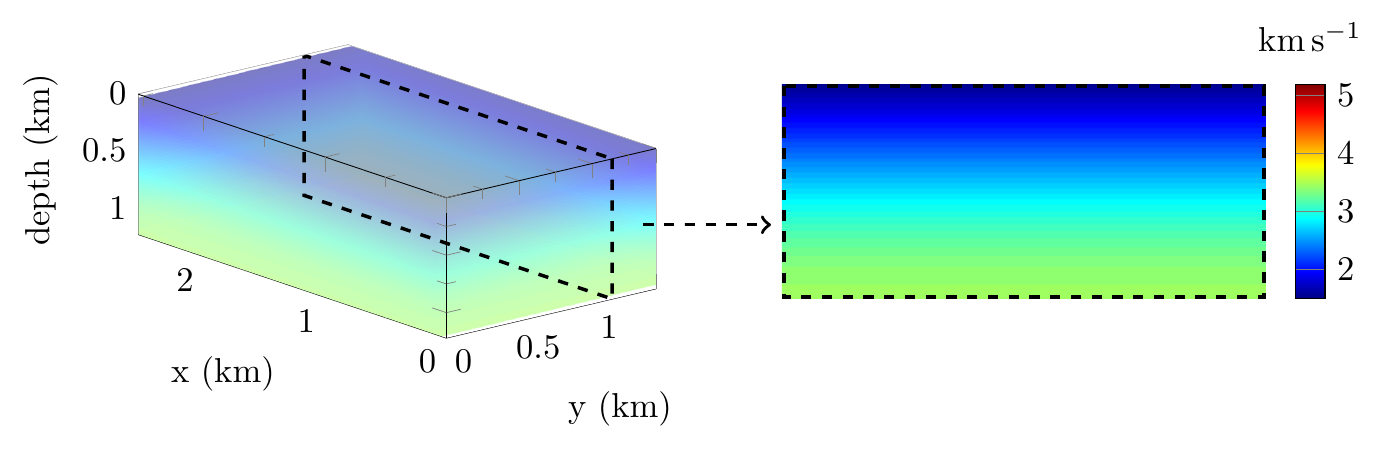}}
\caption{Three-dimensional representation and 
         vertical section at $y=670$\si{\meter} of
         the initial model, which varies in depth only.}
\label{fig:salt_start_1d}
\end{figure}

Time-domain Cauchy data are obtained from this configuration with
$1000$ fixed receivers for each of the $96$ sources. Here the 
sources are positioned on a plane, according to 
Figure~\ref{fig:setup_sketch_c}. The devices are located just underneath 
the surface, at a depth of $10$\si{\meter} for the source and 
$100$\si{\meter} for the receivers. 
We incorporate noise in the time-domain data, with a signal-to-noise
ratio of $15$\si{\decibel}, before we proceed to the Fourier
transform to get the frequency-domain data. 
For the reconstruction we start with an initial model which only varies 
with depth, see Figure~\ref{fig:salt_start_1d}. We do not
assume any contrasting objects in our initial guess, nor 
do we know the value for the background. For the reconstruction 
we only assume the knowledge of the uppermost water layer 
(up to $150$\si{\meter} depth).

To differentiate the simulation set of sources from the 
observation, we reduce their number, change their position
and modify the source wavelet, see 
Table~\ref{table:setup_perspectives_salt}. We perform $100$
iterations of the reconstruction Algorithm~\ref{algo:fwi}, 
with single frequency data at $4$\si{\Hz}. The model 
representation is fixed with $N=1280$ sub-domains where 
piecewise linear functions are employ, for a total of $5120$ 
coefficients. 

\begin{table}[h!]
\begin{center}
\begin{tabular}{|>{\centering\arraybackslash}p{.25\linewidth}|
                 >{\centering\arraybackslash}p{.25\linewidth}|
                 >{\centering\arraybackslash}p{.25\linewidth}|}
 \hline
       & Setup for measurements
       & Setup for simulations        \\ \hline\hline
 Number of sources   &  96            & 60              \\ \hline
 Depth of the sources&  10\si{\meter} & 20\si{\meter} \\ \hline
\end{tabular} \end{center}
\caption{Comparison of acquisition setups employed 
         for the observations and simulations. The 
         source wavelet also differs. For the reconstruction,
         single frequency data at $4$\si{\Hz} and
         a fixed model partition of $N=1280$ are used.}
\label{table:setup_perspectives_salt}
\end{table}

\begin{figure}[h!] 
\centering
{\includegraphics[scale=1]{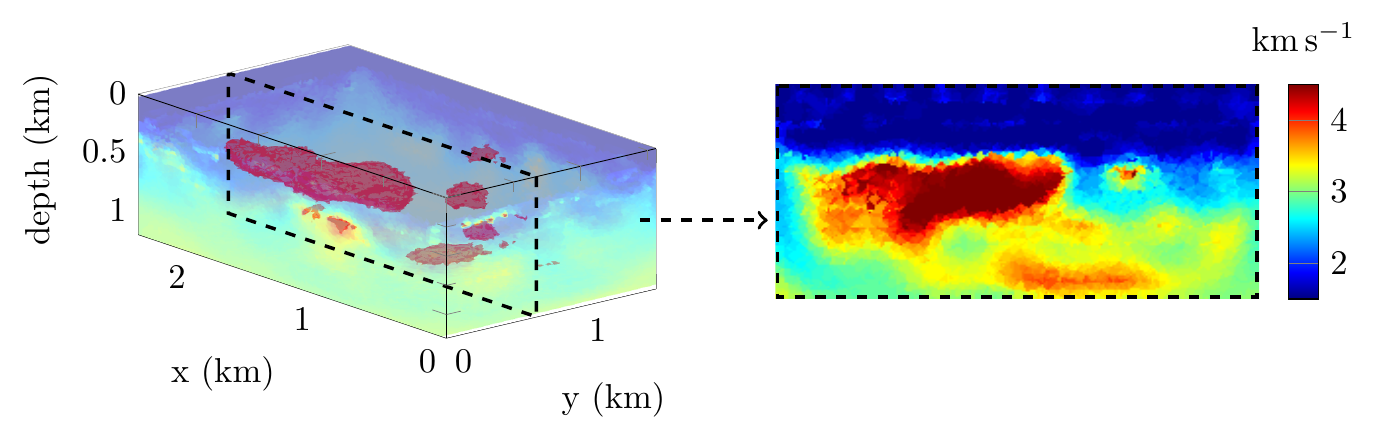}}
\caption{Three-dimensional representation and vertical 
         section at $y=670$\si{\meter} of the 
         reconstruction using $4$\si{\Hz} Cauchy data
         after $100$ iterations and 
         where the simulation setup differs from the original 
         measurement acquisition, see 
         Table~\ref{table:setup_perspectives_salt}.
         The relative $L^2$ difference 
         with the reference model Figure~\ref{fig:salt} 
         is $\mathcal{E}_\text{rel}=0.174$.
         }
\label{fig:fwi_salt_uni-frequency}
\end{figure}

\begin{figure}[h!] 
\centering
{\includegraphics[scale=1]{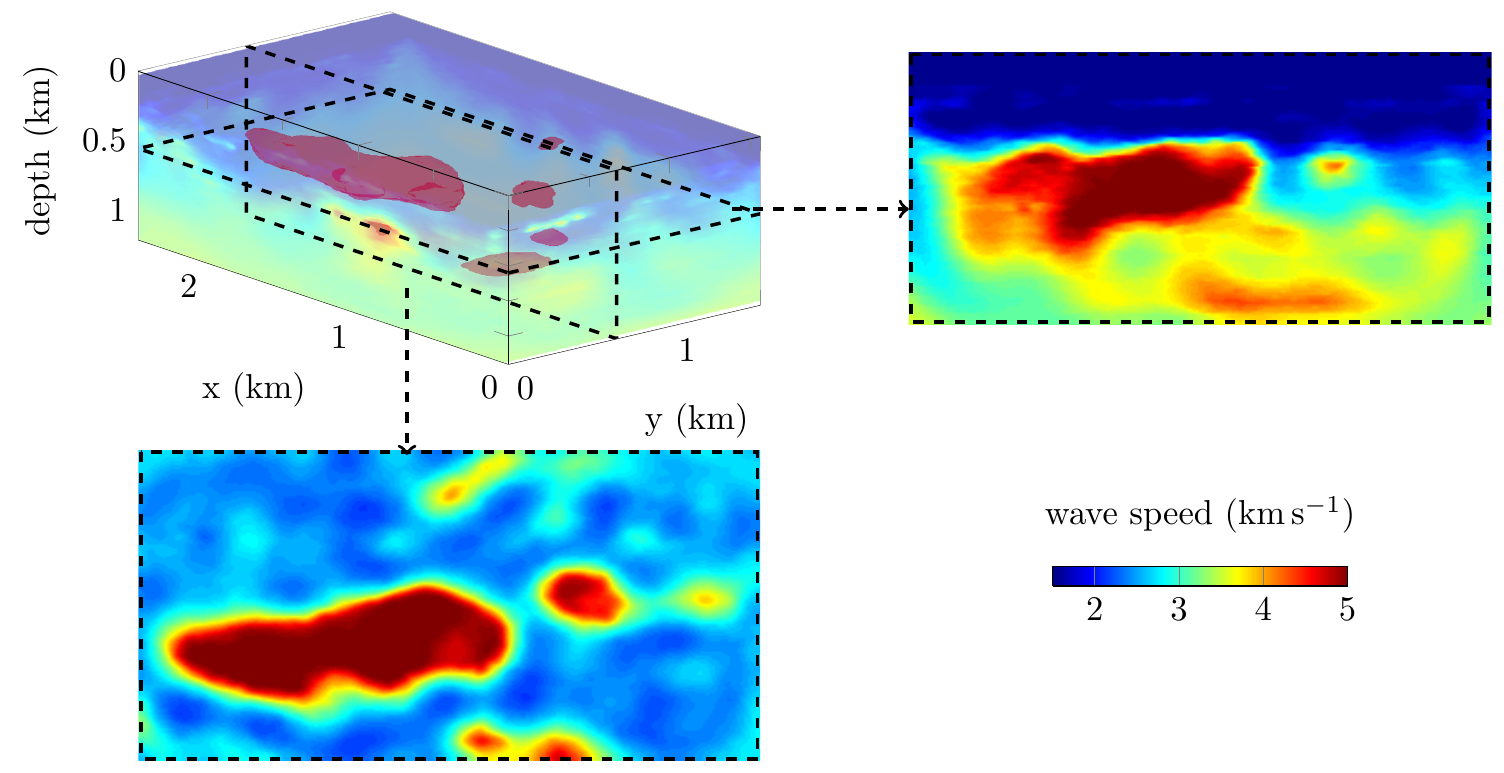}}
\caption{Gaussian filtering (see Remark~\ref{rk:gaussian_image})
         of the reconstruction (Figure~\ref{fig:fwi_salt_uni-frequency})
         obtained using $100$ iterations with $4$\si{\Hz} Cauchy data
         and where the setup for simulation differs from the 
         original measurement acquisition.
         Three-dimensional representation, 
         horizontal section at $550$\si{\meter} 
         depth and vertical section at $y=670$\si{\meter}.
         The relative $L^2$ difference 
         with the reference model Figure~\ref{fig:salt} 
         is $\mathcal{E}_\text{rel}=0.172$.        
         }
\label{fig:fwi_salt_uni-frequency_smooth}
\end{figure}

The reconstruction after $100$ iterations is shown
Figures~\ref{fig:fwi_salt_uni-frequency} 
and~\ref{fig:fwi_salt_uni-frequency_smooth}. 
Despite our initial guess having no information on
the objects, the 
procedure is able to recover the main dome with the 
accurate value, and the shape of smaller domes
(see the horizontal section in 
Figure~\ref{fig:fwi_salt_uni-frequency_smooth}). 
The near boundary information seems missing, as well
as the deepest model variation, due to limited illumination.
However, as we started
with a one dimensional guess and used single frequency
data, the reconstruction is very satisfactory. Once
again, the restricted number of unknowns does not prevent
a good resolution.

We have changed the number of sources in the simulation 
compared to the observation, reducing the numerical cost 
accordingly, yet we make full use of the observed data from the
benefit of our misfit functional defined for Cauchy data.
Furthermore we do not need to know the position of the 
sources employed for observation, nor the source 
wavelet.
The perspective of differentiating the observations and 
simulations acquisition sets is a promising application, 
and appears consistent with the results of 
this preliminary experiment. It would allow less prior
on the observational environment, increasing the robustness
of the procedure, without impacting the resolution of the 
reconstruction.

\newpage
\section{Acknowledgment}

The research of G. Alessandrini and E. Sincich for the 
preparation of this paper has been supported by FRA 2016 
``Problemi inversi, dalla stabilit\`a alla ricostruzione'',
funded by Universit\`a degli Studi di Trieste.
Maarten V. de Hoop acknowledges the Simons Foundation
under the MATH+X program for financial support. He was 
also partially supported by NSF under grant DMS--1559587.
R. Gaburro wishes to acknowledge the support of MACSI, 
the Mathematics Applications Consortium for Science 
and Industry (\url{www.macsi.ul.ie}), funded by the 
Science Foundation Ireland Investigator Award 12/IA/1683.
E. Sincich has also been supported by Gruppo Nazionale per 
l'Analisi Matematica, la Probabilit\`a e le loro Applicazioni 
(GNAMPA) by the grant ``Analisi di problemi inversi: 
stabilit\`a e ricostruzione''.
The research of F. Faucher is supported by the Inria--TOTAL 
strategic action DIP.


\end{document}